\documentclass[12pt]{article}
\usepackage{amsmath}
\usepackage{amssymb}
\renewcommand\smallskip{\vskip\smallskipamount}
\renewcommand\medskip{\vskip\medskipamount}
\renewcommand\bigskip{\vskip\bigskipamount}

\begin{document}

\footnotetext{The first author is partially supported by NSF Grant
DMS-0706794.  The second author is partially supported by NSF
Grant DMS-0707086 and a Sloan Research Fellowship.}

\begin{center}
\begin{large}
\textbf{A Jang Equation Approach to the Penrose Inequality}
\end{large}

\bigskip\smallskip

HUBERT L. BRAY \text{ } \& \text{ } MARCUS A. KHURI

\bigskip
\end{center}
\small

\noindent\textbf{Abstract.}  We introduce a generalized version of
the Jang equation, designed for the general case of the Penrose
Inequality in the setting of an asymptotically flat space-like
hypersurface of a spacetime satisfying the dominat energy
condition.  The appropriate existence and regularity results are
established in the special case of spherically symmetric Cauchy
data, and are applied to give a new proof of the general Penrose
Inequality for these data sets. When appropriately coupled with an
inverse mean curvature flow, analogous existence and regularity
results for the associated system of equations in the nonspherical
setting would yield a proof of the full Penrose Conjecture. Thus
it remains as an important and challenging open problem to
determine whether this system does indeed admit the desired
solutions.
\bigskip

\normalsize

\bigskip

\begin{center}
\textbf{1. Introduction}
\end{center}\setcounter{equation}{0}
\setcounter{section}{1}

\bigskip

  In 1978 P. S. Jang introduced a quasilinear elliptic equation [8],
which Schoen and Yau [13] successfully employed to reduce the
positive mass theorem for general Cauchy data to the case of time
symmetry.  For this reason it has been widely suggested that the
Jang equation could be used in a similar way to reduce the general
Penrose Inequality to the time symmetric case.  However as pointed
out by Malec and \'{O} Murchadha [11], serious issues arise when
one tries to directly apply the steps taken by Schoen and Yau in
[13] (other issues with this process will be pointed out below).
Therefore a new idea is needed, and in this paper is provided in
the form of a generalized Jang equation specifically designed to
treat the Penrose Inequality (PI).\par
   In order to motivate the modification to the Jang equation, let
us recall the precise statement of the PI as well as the suggested
method of proof via the classical Jang equation.  An initial data
set for the Einstein equations is a triple $(M,g,k)$ consisting of
a 3-manifold $M$ (for our purposes with boundary) on which a
positive definite metric $g$ and a symmetric 2-tensor (the
extrinsic curvature) $k$ are defined, which satisfy the constraint
equations
\begin{eqnarray*}
16\pi\mu&=&R-k^{ij}k_{ij}+(g^{ij}k_{ij})^{2},\\
8\pi J_{i}&=&\nabla^{j}(k_{ij}-(g^{ab}k_{ab})g_{ij}),
\end{eqnarray*}
where $R$ is the scalar curvature, $\nabla^{j}$ denotes covariant
differentiation, and $\mu$ and $J_{i}$ are the local matter and
momentum densities respectively. If the initial data are
asymptotically flat, satisfy the dominant energy condition
$\mu\geq |J|_{g}$, and contain an apparent horizon boundary
$\partial M$, then the PI relates the total ADM mass (of a chosen
end) $M_{\mathrm{ADM}}$ to the area $A$ of its outermost minimal
area enclosure (for $\partial M$) by the inequality
\begin{equation}
M_{\mathrm{ADM}}\geq\sqrt{\frac{A}{16\pi}}.
\end{equation}
Furthermore it asserts that if equality holds and the outermost
minimal area enclosure is the boundary of an open bounded domain
$U\subset M$, then $(M-U,g)$ admits an isometric embedding into
the Schwarzschild spacetime with second fundamental form given by
$k$.  The suggested approach for confirming this statement is as
follows.  Look for a surface $\Sigma$ in the product manifold
$(M\times\mathbb{R}, g+dt^{2})$ given by the graph of a function
$t=f(x)$, where $f$ is a solution of Jang's equation.  Then the
induced metric $\overline{g}=g+df^{2}$ on $\Sigma$ has a certain
positivity property for its scalar curvature $\overline{R}$, and
the ADM mass remains unchanged.  One then uses this positivity
property to solve the scalar curvature equation
\begin{equation}
\Delta_{\overline{g}}u-\frac{1}{8}\overline{R}u=0
\end{equation}
on $\Sigma$, to obtain a new metric $u^{4}\overline{g}$ with zero
scalar curvature, and smaller mass.  Thus the hope is that the
area of $\partial\Sigma$ in the new metric $u^{4}\overline{g}$ is
greater than or equal to the area of $\partial M$ in the original
metric $g$, so that an application of the Riemannian PI would give
the desired result.  However, as in [13] it is expected that the
correct boundary behavior for $\Sigma$ is to blow-up and
approximate a cylinder over $\partial M$, but this implies that
the solution $u$ of (1.2) must vanish exponentially fast at
$\partial\Sigma$ (as observed in [11] and [13]) so that we obtain
no contribution from the area of $\partial M$ and hence little
hope of establishing (1.1).  Another failure of this method is
that it has no chance of working in the case of equality, where we
wish to embed the initial data into the Schwarzschild spacetime.
The problem here is that the Schwarzschild spacetime is given by a
warped product metric, and the classical Jang approach only gives
an embedding into a pure product metric.\par
   So we see that there are several problems with the classical
approach to the Penrose Inequality.  The biggest of these problems
is the fact that when the classical Jang surface blows-up inside
the product metric $g+dt^{2}$, it blows-up like a cylinder.  In
other words, the boundary of the Jang surface is infinitely far
away from every point in the surface; this is what causes the
conformal factor to have zero Dirichlet boundary data.  A natural
(and probably first) idea that comes to mind in order to overcome
both difficulties (this one, and the case of equality) is to
consider the warped product metric $g+\phi^{2}dt^{2}$ instead of
the product metric, and require the warping factor $\phi$ to
vanish on $\partial M$. Note that this is compatible with the case
of equality since the warping factor for the Schwarzschild metric
also satisfies this property. We would also like to point out that
although the classical Jang equation has virtually no chance of
establishing the full PI, it has been shown to yield a
Penrose-Like Inequality [9] for general initial data.\par
  Lastly we encourage those who are interested in the current
paper, to compare the motivations and perspective presented here
with the equivalent spacetime formulation presented in [2]. More
precisely, this paper generalizes the Schoen/Yau approach ([13])
to the Positive Mass Theorem in a way which is suitable for the
Penrose Inequality, whereas [2] derives its stimulus from the dual
Lorentzian setting.

\begin{center}
\textbf{2.  The Generalized Jang Equation}
\end{center} \setcounter{equation}{0}
\setcounter{section}{2}

   At this point we see that because of several considerations it
is natural to make the first modification of the Jang approach, by
looking for the Jang surface inside the warped product space
$(M\times\mathbb{R},g+\phi^{2}dt^{2})$.  In order to have any
chance of obtaining a positivity property for the scalar curvature
here, we would like the Jang surface $\Sigma$ to satisfy an
equation with the same structure, namely
\begin{equation}
H_{\Sigma}-\mathrm{Tr}_{\Sigma}K=0,
\end{equation}
where $H_{\Sigma}$ denotes the mean curvature, the tensor $K$ on
$M\times\mathbb{R}$ is an extended version of $k$ from the initial
data, and $\mathrm{Tr}_{\Sigma}K$ denotes the trace of $K$ over
$\Sigma$.  Of course we are free to extend $k$ as we wish. Note
that Schoen and Yau chose to extend $k$ trivially, however as we
will see this extension will not be appropriate for our problem.
The first consideration when looking for a choice of extension, is
that we would like the solutions of Jang's equation to blow-up at
the horizon just as in the classical case, because this gives zero
mean curvature and preserves the area of the horizon inside the
Jang surface. However, it easily seen that the trivial extension
will not allow this in the warped product metric.\par
   Let us consider the extension:
\begin{eqnarray*}
K(\partial_{x^{i}},\partial_{x^{j}})=K(\partial_{x^{j}},\partial_{x^{i}})
&=&k(\partial_{x^{i}},\partial_{x^{j}})\text{ }\text{ }\text{
}\text{ for }\text{ }\text{ }\text{ }1\leq i,j\leq
3,\\
K(\partial_{x^{i}},\partial_{x^{4}})=K(\partial_{x^{4}},\partial_{x^{i}})
&=& 0\!\text{ }\text{ }\text{ }\text{ }\text{ }\text{ }\text{
}\text{ }\text{ }\text{ }\text{ }\text{ }\text{ }\text{ }\text{
}\text{ for }\text{ }\text{ }\text{ }1\leq i\leq
3,\\
K(\partial_{x^{4}},\partial_{x^{4}})&=& k_{44},
\end{eqnarray*}
where $x^{i}$, $i=1,2,3$, are local coordinates on $M$, $x^{4}=t$
is the coordinate on $\mathbb{R}$, and $k_{44}$ is to be
determined.  If the Jang surface blows-up at the horizon
appropriately then it will still approximate a cylinder over the
horizon, but a calculation shows that
\begin{equation*}
H_{\partial M\times\mathbb{R}}=H_{\partial M}+\phi^{-1}\langle
n_{g},\nabla_{g}\phi\rangle_{g_{\phi}} \text{ }\text{ }\text{
}\text{ and }\text{ }\text{ }\text{ }\mathrm{Tr}_{\partial
M\times\mathbb{R}}K=\mathrm{Tr}_{\partial M}k+\phi^{-2}k_{44},
\end{equation*}
where $n_{g}$ is the unit inner normal to $\partial M$ inside
$(M,g)$ and $g_{\phi}$ is the warped product metric. Therefore
since the unit normal to $\Sigma$, given by
\begin{equation*}
N=\frac{\nabla_{g}f-\phi^{-2}\partial_{x^{4}}}{\sqrt{\phi^{-2}+|\nabla_{g}f|^{2}}}
\end{equation*}
where $x^{4}=f(x^{1},x^{2},x^{3})$ expresses $\Sigma$ as a graph,
converges to $\mp n_{g}$ in the process of blowing up to
$\pm\infty$, it is natural to choose
\begin{equation}
k_{44}=\langle N,\phi\nabla_{g}\phi\rangle_{g_{\phi}}
=\frac{\langle\nabla_{g}f,\phi\nabla_{g}\phi\rangle_{g}}{\sqrt{\phi^{-2}+|\nabla_{g}f|^{2}}}
\end{equation}
since $\partial M$ is an apparent horizon and thus satisfies
\begin{equation*}
H_{\partial M}\pm\mathrm{Tr}_{\partial M}k=0.
\end{equation*}
The $\pm$ indicates a future (past) horizon respectively, and the
same expression for $k_{44}$ is valid for both since the Jang
surface will blow-up to $+\infty$ at a future horizon and down to
$-\infty$ at a past horizon.  In other words when $k_{44}$ is
chosen in this way, it is possible for the Jang surface to have
the desired blow-up behavior at horizons.\par
   When the tensor $k$ is extended according to (2.2), we will
refer to equation (2.1) as the generalized Jang equation.  It is
important to note that this particular extension has a natural
interpretation in the dual Lorentzian setting, that is in the
setting of the static spacetime $(M\times\mathbb{R}, g-\phi^{2}
dt^{2})$. More precisely, if we consider the Jang surface $\Sigma$
inside this spacetime then the generalized Jang equation (2.1)
expresses the fact that the second fundamental form of $\Sigma$ in
the Lorentzian setting and the data $k$, when both are pulled back
to the $t=0$ slice, have the same trace over the metric on the
$t=0$ slice (see Appendix B for the relevant calculations).  Thus,
the extension given by (2.2) can be interpreted as the trivial
extension in the dual Lorentzian setting, a fact which is of
paramount importance when proving the rigidity statement in the
case of equality for (1.1).\par
   It turns out that this choice of extension given by (2.2)
actually solves three problems.  Namely, as we have seen it allows
the modified Jang equation to have solutions which blow-up at
horizons, second as we will see later (and eluded to in the
previous paragraph) it is precisely what is needed for the case of
equality, and third it is used to establish a positivity property
for the scalar curvature of the Jang surface $\Sigma$ in the
warped product metric.  The following formula for the Jang surface
$\Sigma$ in the warped product metric is one of the most important
observations of this paper, as it is fundamental for any approach
taken towards the general PI.  As the proof is heavy with
calculation, it is placed in Appendix A.\medskip

\textbf{Theorem 1.}  \textit{Let $\mu$ and $J$ denote the local
energy density and current density associated with the initial
data, respectively.  If the surface $\Sigma$ satisfies the
generalized Jang equation (2.1) and is given by a graph $t=f(x)$,
then its scalar curvature $\overline{R}$ is given by}
\begin{equation}
\overline{R}=16\pi(\mu-J(w))+
|h-K|_{\Sigma}|_{\overline{g}}^{2}+2|q|_{\overline{g}}^{2}
-2\phi^{-1}\mathrm{div}_{\overline{g}}(\phi q),
\end{equation}
\textit{where $h$ is the second fundamental form, $K|_{\Sigma}$ is
the restriction to $\Sigma$ of the extended tensor $K$, $q$ is a
1-form  and $w$ is a vector with $|w|_{g}\leq 1$ given by}
\begin{equation*}
w=\frac{f^{i}\partial_{x^{i}}}{\sqrt{\phi^{-2}+|\nabla_{g}f|^{2}}},\text{
}\text{ }\text{ }\text{ }
q_{i}=\frac{f^{j}}{\sqrt{\phi^{-2}+|\nabla_{g}f|^{2}}}(h_{ij}-(K|_{\Sigma})_{ij}),
\end{equation*}
\textit{with $f^{j}=g^{ij}f_{,i}$.}\medskip

\textbf{Remark.} \textit{The full formula for $\overline{R}$, when
$\Sigma$ does not satisfy any equation, is given in} [2],
\textit{and is referred to as the generalized Schoen-Yau
identity.}\medskip

   Note that (2.3) reduces to the formula obtained by Schoen and
Yau in [13] when $\phi\equiv 1$, which of course corresponds to
the case of the classical Jang equation.  Furthermore the dominant
energy condition ensures that $\mu\geq |J|_{g}$, so that only the
$\mathrm{div}_{\overline{g}}$ term prevents $\overline{R}$ from
being nonnegative.  However, as we shall see, with an appropriate
choice of the warping factor $\phi$ this difficulty can be
overcome to yield the PI, once a full existence theory for the
generalized Jang equation coupled to an inverse mean curvature
flow has been established.

\begin{center}
\textbf{3.  Existence for the Generalized Jang Equation in
Spherical Symmetry}
\end{center} \setcounter{equation}{0}
\setcounter{section}{3}

   In this section we prove the necessary existence and regularity
result needed for the generalized Jang equation, if it is to be
applied to the PI.  We will restrict ourselves to spherically
symmetric initial data.  Therefore the metric $g$ and extrinsic
curvature $k$ have the form
\begin{equation*}
g=g_{11}(r)dr^{2}+\rho^{2}(r)d\Omega^{2},\text{ }\text{ }\text{
}\text{ }k_{ij}=n_{i}n_{j}k_{a}+(g_{ij}-n_{i}n_{j})k_{b},
\end{equation*}
for some functions $g_{11}$, $\rho$, $k_{a}$, $k_{b}$ with the
appropriate fall-off conditions at infinity (to be specified
below), where
\begin{equation*}
n=n^{1}\partial_{r}+n^{2}\partial_{\psi^{2}}+n^{3}\partial_{\psi^{3}}=\sqrt{g^{11}}
\partial_{r}
\end{equation*}
is the unit normal to spheres centered at the origin which will be
denoted by $S_{r}$, and
\begin{equation*}
d\Omega^{2}=(d\psi^{2})^{2}+\sin^{2}\psi^{2}(d\psi^{3})^{2}
\end{equation*}
is the round metric on $\mathbb{S}^{2}$.  We assume that
$M=\mathbb{R}^{3}-B_{0}$ ($B_{0}$ is the ball with boundary
$S_{0}$) so that $\partial M=S_{0}$, with $S_{0}$ an apparent
horizon. Furthermore, we assume that no other apparent horizons
exist in $M$.  This means that the null expansions satisfy
\begin{equation}
\theta_{\pm}=2\left(\sqrt{g^{11}}\frac{\rho_{,r}}{\rho}\pm
k_{b}\right)(r)>0,\text{ }\text{ }\text{ }\text{ }r>0,
\end{equation}
where $\rho_{,r}=\frac{d\rho}{dr}$, and that either
$\theta_{+}(0)=0$, $\theta_{-}(0)=0$, or
$\theta_{+}(0)=\theta_{-}(0)=0$, depending on whether $S_{0}$ is a
future horizon, past horizon, or both, respectively.\par
   We now derive the generalized Jang equation.  Let the Jang surface
$\Sigma$ be given as the graph of a function $t=f(r)$, then the
unit normal to $\Sigma$ in the warped product metric
$g_{\phi}=g+\phi^{2}dt^{2}$ is
\begin{equation*}
N=\frac{g^{11}f_{,r}\partial_{r}-\phi^{-2}\partial_{t}}{\sqrt{\phi^{-2}
+g^{11}f_{,r}^{2}}}:=
N^{1}\partial_{r}+N^{2}\partial_{\psi^{2}}+N^{3}\partial_{\psi^{3}}
+N^{4}\partial_{t}.
\end{equation*}
The mean curvature of $\Sigma$ with respect to $N$ is
\begin{equation*}
H_{\Sigma}=\sum_{i=1}^{4}N_{;i}^{i}=\sum_{i=1}^{4}\frac{1}{\sqrt{|g_{\phi}|}}
\partial_{i}\left(\sqrt{|g_{\phi}|}N^{i}\right)=\frac{\sqrt{g^{11}}}{\phi}
\partial_{r}(\phi v)+2\sqrt{g^{11}}\frac{\rho_{,r}}{\rho}v
\end{equation*}
where $|g_{\phi}|=\det g_{\phi}$, $N_{;i}^{i}$ denotes covariant
differentiation with respect to $g_{\phi}$, and
\begin{equation*}
v=\frac{\phi\sqrt{g^{11}}f_{,r}}{\sqrt{1+\phi^{2}
g^{11}f_{,r}^{2}}}.
\end{equation*}
Furthermore the extension $K$ of the extrinsic curvature given by
(2.2) requires that
\begin{equation*}
k_{44}=\frac{g^{11}\phi\phi_{,r}f_{,r}}{\sqrt{\phi^{-2}+g^{11}
f_{,r}^{2}}},
\end{equation*}
so that if $\overline{g}=g+\phi^{2}df^{2}$ denotes the induced
metric on $\Sigma$ we have
\begin{eqnarray*}
\mathrm{Tr}_{\Sigma}K &=& \overline{g}^{ij}k_{ij}
+\overline{g}^{ij}f_{,i}f_{,j}k_{44}\\
&=&
k_{a}+2k_{b}-\frac{g^{11}f_{,r}^{2}}{\phi^{-2}+g^{11}f_{,r}^{2}}k_{a}
+\left(\frac{g^{11}f_{,r}^{2}}{1+\phi^{2}g^{11}f_{,r}^{2}}\right)
\left(\frac{g^{11}\phi^{2}\phi_{,r}f_{,r}}{\sqrt{1+\phi^{2}
g^{11}f_{,r}^{2}}}\right)\\
&=&(1-v^{2})k_{a}+2k_{b}+\sqrt{g^{11}}\frac{\phi_{,r}}{\phi}
v^{3}.
\end{eqnarray*}
Thus the generalized Jang equation (2.1) takes the form
\begin{equation}
\sqrt{g^{11}}v_{,r}+2\left(\sqrt{g^{11}}\frac{\rho_{,r}}{\rho}v-k_{b}\right)
+(v^{2}-1)k_{a}+\sqrt{g^{11}}v\frac{\phi_{,r}}{\phi}(1-v^{2})=0.
\end{equation}\par
   For the proof of the PI in the next section we will need to set
\begin{equation}
\phi=\rho_{,s}=\frac{\sqrt{1-v^{2}}}{\sqrt{g_{11}}}\rho_{,r}
\end{equation}
where $s$ is the radial arclength parameter in the $\overline{g}$
metric, that is
\begin{equation*}
s=\int_{0}^{r}\sqrt{g_{11}+\phi^{2}
f_{,r}^{2}}=\int_{0}^{r}\frac{\sqrt{1-v^{2}}}{\sqrt{g_{11}}}.
\end{equation*}
Thus our existence results shall only concern the case in which
$\phi$ is given by (3.3).  First note that $\rho_{,r}(r)>0$,
$r>0$, since the condition (3.1) shows that
\begin{equation*}
2\sqrt{g^{11}}\frac{\rho_{,r}}{\rho}(r)=H_{S_{r},g}=\frac{1}{2}(\theta_{+}+\theta_{-})(r)>0
,\text{ }\text{ }\text{ }\text{ }r>0,
\end{equation*}
so $\phi$ is well-defined.  Secondly, when $\phi$ is given by
(3.3) we have
\begin{equation*}
\sqrt{g^{11}}v\frac{\phi_{,r}}{\phi}(1-v^{2})=
v\frac{\phi_{,s}}{\phi}\sqrt{1-v^{2}}
=v\frac{\rho_{,ss}}{\rho_{,s}}\sqrt{1-v^{2}},
\end{equation*}
and
\begin{equation*}
\rho_{,ss}=\rho_{,rr}\left(\frac{1-v^{2}}{g_{11}}\right)
-\frac{vv_{,r}}{g_{11}}\rho_{,r}-\frac{1}{2}\frac{g_{11,r}}{g_{11}^{2}}(1-v^{2})
\rho_{,r}.
\end{equation*}
Therefore, the generalized Jang equation (3.2) becomes
\begin{equation}
\sqrt{g^{11}}(1-v^{2})v_{,r}+(1-v^{2})F_{\mp}(r,v)\pm\theta_{\mp}=0
\end{equation}
where
\begin{equation*}
F_{\mp}(r,v)=\mp 2\sqrt{g^{11}}\frac{\rho_{,r}}{\rho}\frac{1}{1\pm
v}-k_{a}+\frac{v}{\sqrt{g^{11}}}\frac{\rho_{,rr}}{\rho_{,r}}-\frac{v}{2}
\frac{g_{11,r}}{g_{11}^{2}}.
\end{equation*}\par
   It is interesting to note the role of the null expansions in
equation (3.4).  It turns out that the outermost apparent horizon
hypothesis (3.1) is the primary reason that we are able to obtain
an existence and regularity result in $M$.  This is analogous with
the theory developed by Schoen and Yau in [13] for the classical
Jang equation, in that the absence of horizons leads to
regularity.  Furthermore observe that according to the definition
of $v$, $v=\pm 1$ corresponds to blow-up of the Jang surface
$\Sigma$.  Thus for us, by regularity of a solution to (3.4) we
mean not only that the solution possesses a large number of
continuous derivatives, but that it satisfies $-1<v<1$ as well.
The following existence result is what we require for the PI in
the next section.  Assume that the initial data satisfy the
following fall-off conditions as $r\rightarrow\infty$:
\begin{eqnarray}
|k(r)|_{g}\!\!\!&\leq&\!\!\! Cr^{-2},\text{ }\text{ }\text{
}\text{ }
|\mathrm{Tr}_{g}k(r)|\leq Cr^{-3},\\
|(g_{11}-1)(r)|+r|g_{11,r}(r)|\!\!\!&\leq&\!\!\! Cr^{-1},\text{
}\text{ }\text{ }\text{ }
|\rho(r)-r|+r|\rho_{,r}(r)-1|+r^{2}|\rho_{,rr}(r)|\leq C,\nonumber
\end{eqnarray}
for some constant $C$.\medskip

\textbf{Theorem 2.}  \textit{Assume that the initial data are
smooth, satisfy the outermost apparent horizon condition (3.1),
and the asymptotics (3.5).  Then given $\alpha\in(-1,1)$ there
exists a unique solution $v\in C^{\infty}((0,\infty))\cap
C^{1}([0,\infty))$ of (3.4) such that $-1<v(r)<1$, $r>0$, and
$v(0)=\alpha$.  If $S_{0}$ is a past (future) horizon then the
same conclusion holds with $v(0)=\pm 1$, respectively. Furthermore
the solution $v$ has the following asymptotics}
\begin{equation*}
|v(r)|+r|v_{,r}(r)|\leq Cr^{-2},\text{ }\text{ }\textit{ as
}\text{ }\text{ }r\rightarrow\infty,
\end{equation*}
\textit{for a constant $C$ depending only on
$|g|_{C^{1}((0,\infty))}$ and $|k|_{C^{0}((0,\infty))}$.}\medskip

\textbf{Remark.}  \textit{The generalized Jang equation addressed
here corresponds to a specific choice of $\phi$, namely that given
by (3.3).}\medskip

\textit{Proof.}  We first establish the fundamental a priori
estimate
\begin{equation}
-1<v(r)<1,\text{ }\text{ }\text{ }\text{ }r>0,
\end{equation}
as a consequence of the outermost apparent horizon condition
(3.1).  First consider the case when $|v(0)|<1$.  Then arguing by
contradiction there must exist a smallest value $r_{0}>0$ such
that $v(r_{0})=1$.  It follows that there is an $\varepsilon>0$
such that
\begin{equation*}
v(r)<1,\text{ }\text{ }\text{ }\text{ }v_{,r}(r)\geq 0,\text{
}\text{ }\text{ }\text{ }r_{0}-r<\varepsilon.
\end{equation*}
However from equation (3.4) and (3.1) we have
\begin{equation*}
\sqrt{g^{11}}(1-v^{2})v_{,r}(\overline{r})=-\theta_{-}(\overline{r})
-(1-v^{2})F_{-}(\overline{r},v)<0\text{ }\text{ }\text{ for some
}\text{ }\text{ }r_{0}-\overline{r}<\varepsilon.
\end{equation*}
A similar argument can be used if $v(r_{0})=-1$ by replacing
$\theta_{-}$ with $\theta_{+}$.  This establishes (3.6) if
$|v(0)|<1$.  If $v(0)=1$ and $S_{0}$ is a past horizon then these
same arguments yield (3.6) as long as $v_{,r}(0)<0$ (similarly if
$v(0)=-1$ and $S_{0}$ is a future horizon then we need
$v_{,r}(0)>0$).  To confirm this it suffices to write out partial
Taylor expansions at $r=0$ for all functions appearing in (3.4).
It follows that
\begin{equation}
\sqrt{g^{11}}(0)v_{,r}(0)^{2}+F_{-}(0,v)v_{,r}(0)-\frac{1}{2}\theta_{-,r}(0)=0.
\end{equation}
We can assume without loss of generality that $\theta_{-,r}(0)>0$
(which is of course consistent with (3.1)), by slightly perturbing
the initial data.  Therefore we find that $v_{,r}(0)<0$ as
desired.  This establishes (3.6).\par
   In order to prove existence, we need to obtain a priori
estimates for the derivatives of $v$.  The first task in this
direction is to improve (3.6).  Let $(r_{0},r_{1})$ be an interval
on which $v(r)\geq 0$ and $v_{,r}(r)\leq 0$, then
\begin{equation}
v(r)\leq v(r_{0}),\text{ }\text{ }\text{ }\text{ }r\in
(r_{0},r_{1}).
\end{equation}
If $(r_{0},r_{1})$ is an interval on which $v(r)\geq 0$ and
$v_{,r}(r)\geq 0$ then from equation (3.4) and (3.6) we have
\begin{equation*}
0\leq -\theta_{-}(r)+\overline{F}_{-}(r)(1-v^{2}),\text{ }\text{
}\text{ }\text{ }r\in (r_{0},r_{1}),
\end{equation*}
where
\begin{equation*}
|F_{-}(r)|\leq \overline{F}_{-}(r)\text{ }\text{ }\text{ with
}\text{ }\text{ }C\leq\overline{F}_{-}(r)\text{ }\text{ }\text{
and }\text{ }\text{
}\frac{\theta_{-}(r)}{\overline{F}_{-}(r)}\rightarrow C^{-1}\text{
}\text{ }\text{ as }\text{ }\text{ }r\rightarrow\infty,
\end{equation*}
for some universal constant $C>0$.  Then according to (3.1) and
the fall-off conditions (3.5),
\begin{equation}
(1-v^{2})\geq\frac{\theta_{-}(r)}{\overline{F}_{-}(r)}\geq\delta>0
\end{equation}
for some $0<\delta<1$ independent of the interval $(r_{0},r_{1})$
if $r_{0}\geq\varepsilon>0$.  Similar estimates can be obtained
when $v(r)\leq 0$. Thus by combining (3.8) and (3.9) we conclude
that there exists $0<\delta(\alpha)<1$ for each $|\alpha=v(0)|<1$
such that
\begin{equation}
-1+\delta(\alpha)\leq v(r)\leq 1-\delta(\alpha),\text{ }\text{
}\text{ }\text{ }r\in [0,\infty).
\end{equation}
If $v(0)=\pm 1$ and $S_{0}$ is a past (future) horizon then
$v_{,r}(0)<(>) 0$ from (3.7), so the same arguments provide
$0<\delta(\varepsilon)<1$ for each $\varepsilon>0$ such that
\begin{equation}
-1+\delta(\varepsilon)\leq v(r)\leq 1-\delta(\varepsilon),\text{
}\text{ }\text{ }\text{ }r\in [\varepsilon,\infty).
\end{equation}\par
   With the aid of (3.10) and (3.11) we can now simply
differentiate equation (3.4) to inductively show that there exist
constants $C(l,\alpha)$, $|\alpha=v(0)|<1$, $l\in\mathbb{Z}_{+}$
such that
\begin{equation*}
|v|_{C^{l}([0,\infty))}\leq C(l,\alpha),
\end{equation*}
and constants $C(l,\varepsilon)$, $\varepsilon>0$, if $v(0)=\pm 1$
and $S_{0}$ is a past (future) horizon, such that
\begin{equation*}
|v|_{C^{l}([\varepsilon,\infty))}\leq C(l,\varepsilon).
\end{equation*}
Moreover because we can solve for $v_{,r}(0)$ from (3.7), we also
obtain the global $C^{1}$ estimate
\begin{equation*}
|v|_{C^{1}([0,\infty))}\leq C.
\end{equation*}
At this point we can then make a standard application of the
Leray-Schauder fixed point theorem (or alternatively the method of
continuity) to obtain a global solution $v\in
C^{\infty}((0,\infty))\cap C^{1}([0,\infty))$ with prescribed
$v(0)\in [-1,1]$.  Of course if $v(0)=\pm 1$ then we require
$S_{0}$ to be a past (future) horizon respectively.\par
   Lastly we show that $v$ has the correct asymptotics at
infinity.  By (3.10), (3.11), and the fall-off conditions (3.5) we
can write equation (3.4) as
\begin{equation}
v_{,r}+\frac{2r^{-1}}{1-v^{2}}v=O(r^{-3}+r^{-2}v),\text{ }\text{
}\text{ }\text{ }0<r_{0}<r<\infty,
\end{equation}
noting that
\begin{equation*}
\mathrm{Tr}_{g}k=k_{a}+2k_{b}.
\end{equation*}
Thus the solution on the interval $(r_{0},\infty)$ can be
represented by
\begin{equation*}
v(r)=\exp\left(-\int_{r_{0}}^{r}\frac{2r^{-1}}{1-v^{2}}\right)
\left[\int_{r_{0}}^{r}O(r^{-3}+r^{-2}v)\exp\left(\int_{r_{0}}^{r}\frac{2r^{-1}}{1-v^{2}}\right)
+v(r_{0})\right].
\end{equation*}
It follows that
\begin{equation*}
|v(r)|\leq Cr^{-1}\text{ }\text{ }\text{ as }\text{ }\text{
}r\rightarrow\infty.
\end{equation*}
Plugging this back into the above representation produces
\begin{equation}
|v(r)|\leq Cr^{-2}\text{ }\text{ }\text{ as }\text{ }\text{
}r\rightarrow\infty.
\end{equation}
Therefore from (3.12) and (3.13) we have
\begin{equation*}
|v_{,r}(r)|\leq C(r^{-3}+r^{-1}|v(r)|)\leq Cr^{-3}\text{ }\text{
}\text{ as }\text{ }\text{ }r\rightarrow\infty.
\end{equation*}
Q.E.D.\par

\begin{center}
\textbf{4.  Proof of the Penrose Inequality in the Case of
Spherical Symmetry}
\end{center} \setcounter{equation}{0}
\setcounter{section}{4}

   In this section we show how to apply the generalized Jang equation
to treat the PI as stated in section $\S 1$.  The proof presented
here is restricted to the case of spherically symmetric initial
data. However as illustrated in the next section, this method
could be generalized to cover arbitrary data if an analogous
existence result for the generalized Jang equation is established.
A significant difference in the general case is that it is
necessary to solve a system of equations (see [2]), whereas in the
case of spherical symmetry only the generalized Jang equation need
be solved, as the system actually decouples. Note that several
different proofs of the PI for spherically symmetric initial data
have been put forward (eg. [4], [5], [7], [10]).  The proof
presented below is new, and appears to be novel in that it has the
potential to generalize.\par
   Using notation already established in the previous section, let
$\overline{g}=g+\phi^{2}df^{2}$ denote the induced metric on the
Jang surface and write
\begin{equation*}
\overline{g}=ds^{2}+\rho^{2}(s)d\Omega^{2},
\end{equation*}
where
\begin{equation*}
s=\int_{0}^{r}\sqrt{g_{11}+\phi^{2}
f_{,r}^{2}}=\int_{0}^{r}\frac{\sqrt{1-v^{2}}}{\sqrt{g_{11}}}
\end{equation*}
is the radial arclength parameter in the $\overline{g}$ metric. We
first derive the Hawking mass.  By trying to transform
$\overline{g}$ into a Schwarzschild metric we have
\begin{equation*}
\overline{g}=\frac{1}{\rho_{,s}^{2}}d\rho^{2}+\rho^{2}d\Omega^{2}
=\left(1-\frac{2m(s)}{\rho}\right)^{-1}d\rho^{2}+\rho^{2}d\Omega^{2}
\end{equation*}
for some function $m(s)$.  Solving for $m(s)$ produces
\begin{equation*}
2m=\rho(1-\rho_{,s}^{2}),
\end{equation*}
where
\begin{equation*}
\rho_{,s}=\frac{\sqrt{1-v^{2}}}{\sqrt{g_{11}}}\rho_{,r}.
\end{equation*}
As in the previous section let $S_{r}$ denote a sphere of radius
$r$, then
\begin{equation}
A_{\overline{g}}(S_{r})=A_{g}(S_{r})=4\pi\rho^{2},\text{ }\text{
}\text{ }\text{
}H_{S_{r},\overline{g}}=2\frac{\rho_{,s}}{\rho}=2\frac{\sqrt{1-v^{2}}}{\sqrt{g_{11}}}
\frac{\rho_{,r}}{\rho},
\end{equation}
where $A_{\overline{g}}(S_{r})$ and $H_{S_{r},\overline{g}}$
denote area and mean curvature in the $\overline{g}$ metric
respectively.  It follows that
\begin{equation*}
m(r)=\frac{1}{2}\rho(1-\rho_{,s}^{2})
=\sqrt{\frac{A_{\overline{g}}(S_{r})}{16\pi}}\left(1-\frac{1}{16\pi}
\int_{S_{r}}H^{2}_{S_{r},\overline{g}}d\sigma_{\overline{g}}\right)
\end{equation*}
is precisely the Hawking mass with $d\sigma_{\overline{g}}$
representing the area form in the $\overline{g}$ metric.
Furthermore a calculation shows that the scalar curvature of
$\overline{g}$ is given by
\begin{equation*}
\overline{R}=2\rho^{-2}(1-2\rho\rho_{ss}-\rho_{s}^{2}),
\end{equation*}
and therefore
\begin{equation}
2m_{,s}=\rho_{,s}-\rho_{,s}^{3}-2\rho\rho_{,ss}\rho_{,s}=\frac{1}{2}\rho_{,s}\rho^{2}\overline{R}.
\end{equation}\par
   We will now use these formulas to obtain the PI.  Set the
warping factor by $\phi=\rho_{,s}$ so that Theorem 2 guarantees a
unique solution of the generalized Jang equation.  We also assume
that $S_{0}$ is a past horizon so that we can take $v(0)=1$ in
Theorem 2 (the same arguments below will work if $S_{0}$ is a
future horizon).  Note that the outermost apparent horizon
condition (3.1) guarantees that
\begin{equation*}
H_{S_{r},g}=\frac{1}{2}(\theta_{+}+\theta_{-})(r)>0,\text{ }\text{
}\text{ }\text{ }r>0,
\end{equation*}
which implies
\begin{equation}
\phi(r)=\rho_{,s}(r)=\frac{\sqrt{1-v^{2}}}{\sqrt{g_{11}}}\rho_{,r}(r)
=\frac{\sqrt{1-v^{2}}}{2}\rho H_{S_{r},g}>0,\text{ }\text{
}r>0,\text{ }\text{ }\phi(0)=0,
\end{equation}
where we have also used the estimate (3.6) and $v(0)=1$.  Now in
order to obtain the PI just integrate equation (4.2):
\begin{equation*}
m(\infty)-m(0)=\int_{0}^{\infty}m_{,s}ds=\int_{0}^{\infty}\rho_{,s}\rho^{2}\overline{R}
ds=\frac{1}{4\pi}\int_{\Sigma}\rho_{,s}\overline{R}d\omega_{\overline{g}}
\end{equation*}
by (4.1), where $d\omega_{\overline{g}}$ is the volume form on the
Jang surface $\Sigma$.  Then applying the formula for
$\overline{R}$ from Theorem 1, the dominant energy condition, the
definition of $\phi$, as well as the divergence theorem and (4.3),
we have
\begin{eqnarray}
m(\infty)-m(0)&=&\frac{1}{4\pi}\int_{\Sigma}\rho_{,s}(2(\mu-J(w))+
|h-K|_{\Sigma}|_{\overline{g}}^{2}+2|q|_{\overline{g}}^{2})d\omega_{\overline{g}}\nonumber\\
&
&-\frac{1}{2\pi}\int_{\Sigma}\mathrm{div}_{\overline{g}}(\phi q)d\omega_{\overline{g}}\\
&\geq&-\frac{1}{2\pi}\int_{\partial\Sigma\cup\partial\infty}
\phi\overline{g}(q,n_{\overline{g}})d\sigma_{\overline{g}},\nonumber
\end{eqnarray}
where $K|_{\Sigma}$ is the restriction to $\Sigma$ of the extended
(by (2.2)) tensor $K$, $n_{\overline{g}}$ is the unit outer normal
(as a 1-form), and $q$, $w$ are given in Theorem 1.  According to
a calculation relegated to Appendix C
\begin{equation*}
\phi\overline{g}(q,n_{\overline{g}})d\sigma_{\overline{g}}
=\pm\frac{2\rho_{,r}v}{\sqrt{g_{11}}}
\left(\sqrt{g^{11}}\frac{\rho_{,r}}{\rho}v-k_{b}\right)\rho^{2}d\sigma,
\end{equation*}
where $d\sigma$ is the Euclidean area element.  Therefore the
boundary integral of (4.4) taken over $\partial\Sigma$ is zero,
since $v(0)=1$ and $S_{0}$ is a past horizon.  Also the boundary
integral over $\partial\infty$ is zero as well according to the
asymptotics for $v(r)$ given in Theorem 2 and the fall-off
conditions (3.5).  It follows that
\begin{equation*}
M_{\mathrm{ADM}}-\sqrt{\frac{A_{g}(S_{r})}{16\pi}}=m(\infty)-m(0)\geq
0,
\end{equation*}
since (4.1) and $v(0)=1$ give $H_{S_{0},\overline{g}}=0$.\par
   Lastly we prove the rigidity statement in the case of equality.
The same arguments can be used to deal with the
$\mathrm{div}_{\overline{g}}$ term in (4.4), thus
\begin{equation*}
0=M_{\mathrm{ADM}}-\sqrt{\frac{A_{g}(S_{r})}{16\pi}}\geq\frac{1}{4\pi}\int_{\Sigma}
\rho_{,s}(2(\mu-|J|_{g})+ |h-K|_{\Sigma}|_{\overline{g}}^{2}
+2|q|_{\overline{g}}^{2})d\omega_{\overline{g}}.
\end{equation*}
Hence
\begin{equation*}
\mu-|J|_{g}\equiv 0,\text{ }\text{ }\text{ }\text{
}h-K|_{\Sigma}\equiv 0,\text{ }\text{ }\text{ }\text{ }q\equiv 0.
\end{equation*}
It then follows from Theorem 1 that $\overline{R}\equiv 0$.  We
can now apply the time symmetric PI to the Jang surface $\Sigma$
to obtain $\overline{g}\cong g_{\mathrm{SC}}$, that is
$\overline{g}$ is isometric to the standard slice of the
Schwarzschild spacetime
\begin{equation*}
g_{\mathrm{SC}}=\left(1-\frac{2M_{\mathrm{ADM}}}{r}\right)^{-1}dr^{2}+r^{2}d\Omega^{2}.
\end{equation*}
Hence
\begin{equation*}
\rho=r,\text{ }\text{ }\text{ }\text{
}\overline{g}_{11}=\left(1-\frac{2M_{\mathrm{ADM}}}{r}\right)^{-1},
\end{equation*}
so that
\begin{equation*}
\phi=\rho_{,s}=\frac{1}{\sqrt{\overline{g}_{11}}}\rho_{,r}
=\left(1-\frac{2M_{\mathrm{ADM}}}{r}\right)^{1/2}.
\end{equation*}
This says that $\phi$ is the correct warping factor for the
Schwarzschild spacetime, and furthermore since
\begin{equation*}
g=\overline{g}-\phi^{2}df^{2}=g_{\mathrm{SC}}-\phi^{2}df^{2},
\end{equation*}
the graph map $G:M\rightarrow\mathbb{SC}^{4}$ provides an
isometric embedding of the initial data $(M,g)$ into the
Schwarzschild spacetime $(\mathbb{SC}^{4},g_{\mathrm{SC}}-\phi^{2}
dt^{2})$.  Finally a calculation (Appendix B) shows that
$h-K|_{\Sigma}\equiv 0$ implies that the second fundamental form
of $G(M)\subset\mathbb{SC}^{4}$ is precisely given by the initial
data $k$.

\pagebreak

\begin{center}
\textbf{5.  Approach to the General Case}
\end{center} \setcounter{equation}{0}
\setcounter{section}{5}

Here we discuss how the approach of the previous section may be
generalized to the case of arbitrary initial data (without
spherical symmetry), whenever appropriate solutions exist to a
canonical system of equations constructed from the generalized
Jang equation and the inverse mean curvature flow. Whether or not
such solutions do indeed exist is thus a very important and
challenging open problem. To proceed, we assume without loss of
generality that $\partial M$ is an outermost apparent horizon.
More precisely, none of the components of $\partial M$ are
separated from spatial infinity by another horizon. The idea is
that if $\partial M$ is not outermost, then we should replace $M$
with the submanifold $\widetilde{M}$ such that
$\partial\widetilde{M}$ is an outermost horizon.  Then the PI for
$\widetilde{M}$ implies the PI for $M$. This assumption is made in
order to facilitate the existence of a smooth Jang surface
$\Sigma$ on the interior of $M$, which also blows-up at the
boundary (this boundary behavior guarantees that $\partial\Sigma$
is minimal). Such solutions have been shown to exist in [12], at
least for the classical Jang equation. Furthermore, we assume the
existence of a smooth Inverse Mean Curvature Flow (IMCF) inside
the Jang surface, starting from any one of the components of the
outermost minimal surface enclosing $\partial M$ (it is customary
to take the one with largest area). As the name suggests, IMCF
refers to the flow of 2-surfaces in $\Sigma$ in which the surfaces
flow in the outward normal direction at a rate equal to the
inverse of their mean curvatures at each point. Originally
introduced by Geroch [3], this flow has been generalized and used
successfully by Huisken and Ilmanen [6] to prove the PI in the
time symmetric case. Note that since $\partial\Sigma$ is minimal,
the existence of an outermost minimal surface is guaranteed ([6],
Lemma 4.1), and moreover the region between the outermost minimal
surface and spatial infinity, denoted by $\widetilde{\Sigma}$,
contains no other compact minimal surfaces, and each component of
$\widetilde{\Sigma}$ has spherical topology.  This observation is
required so that the weak formulation of IMCF given by Huisken and
Ilmanen has a smooth start at the boundary.  Although use of the
weak formulation is necessary, since it is not difficult to find
examples where the flow develops singularities, for the sake of
simplicity of exposition, in this paper we assume the existence of
a smooth flow.  This means that in $\widetilde{\Sigma}$, the
induced metric $\overline{g}$ may be written as
\begin{equation*}
\overline{g}=H_{S_{r},\overline{g}}^{-2}dr^{2}+\sum_{i,j=1}^{2}\widehat{g}_{ij}d\theta^{i}d\theta^{j},
\end{equation*}
where the surfaces $r=\mathrm{const.}$ are the flow surfaces
denoted by $S_{r}$ (each having spherical topology), and
$\theta^{i}$ are local coordinates on $S_{r}$.  We also set
$S_{0}$ to be a component of the outermost minimal surface
$\partial\widetilde{\Sigma}$.\par
   Let $m(r)$ again be the Hawking mass, then a
well-known formula (due to Geroch [3]) gives
\begin{equation*}
\frac{dm}{dr}(r)=\sqrt{\frac{A_{\overline{g}}(S_{r})}{16\pi}}\left[\frac{1}{2}
+\frac{1}{16\pi}\int_{S_{r}}\left(2\frac{|\nabla_{S_{r}}H_{S_{r},\overline{g}}|^{2}}
{H_{S_{r},\overline{g}}^{2}}+\overline{R}-2K_{S_{r}}+\frac{1}{2}(\lambda_{1}-\lambda_{2})^{2}
\right)d\sigma_{\overline{g}}\right]
\end{equation*}
where $K_{S_{r}}$ is the Gaussian curvature of $S_{r}$ and
$\lambda_{1}$, $\lambda_{2}$ are its principle curvatures.  Since
each $S_{r}$ has spherical topology (as the flow is assumed to be
smooth), the Gauss-Bonnet Theorem shows that
\begin{equation*}
\int_{S_{r}}K_{S_{r}}d\sigma_{\overline{g}}=4\pi.
\end{equation*}
We then have
\begin{eqnarray}
m(\infty)-m(0)&=&\int_{0}^{\infty}\frac{dm}{dr}dr\\
&\geq&\frac{1}{(16\pi)^{3/2}}\int_{0}^{\infty}\left( \int_{S_{r}}
\sqrt{A_{\overline{g}}(S_{r})}\overline{R}d\sigma_{\overline{g}}\right)dr\nonumber\\
&=&\frac{1}{(16\pi)^{3/2}}\int_{\widetilde{\Sigma}}\sqrt{A_{\overline{g}}(S_{r})}
H_{S_{r},\overline{g}}\overline{R}d\omega_{\overline{g}}\nonumber\\
&\geq&-\frac{2}{(16\pi)^{3/2}}\int_{\widetilde{\Sigma}}\frac{\sqrt{A_{\overline{g}}(S_{r})}
H_{S_{r},\overline{g}}}{\phi}\mathrm{div}_{\overline{g}}(\phi q)
d\omega_{\overline{g}}\nonumber
\end{eqnarray}
according to (2.3), the coarea formula, and the fact that
$H_{S_{r},\overline{g}}>0$ under smooth IMCF.  This motivates the
choice
\begin{equation}
\phi=\sqrt{A_{\overline{g}}(S_{r})}H_{S_{r},\overline{g}},
\end{equation}
since an application of the divergence theorem then yields
\begin{equation}
m(\infty)-m(0)\geq-\frac{2}{(16\pi)^{3/2}}\int_{S_{0}\cup\partial\infty}
\phi\overline{g}(q,n_{\overline{g}})d\sigma_{\overline{g}},
\end{equation}
where $n_{\overline{g}}$ is the unit outer normal with respect to
$\overline{g}$.\par
  In order to obtain the PI from (5.3), we note that since the
solution of the generalized Jang equation vanishes very fast at
spatial infinity, $m(\infty)$ is the original ADM mass of $M$ and
the integral at $\partial\infty$ vanishes (as $\phi$ remains
bounded). Furthermore because $S_{0}$ is a minimal surface
$m(0)=\sqrt{A_{\overline{g}}(S_{0})/16\pi}$.  If $S_{0}$ does not
intersect $\partial\Sigma$ then the solution of the generalized
Jang equation remains bounded on $S_{0}$, so $\phi|_{S_{0}}=0$
implies that the boundary integral vanishes in this case.  On the
other hand, if a portion of $S_{0}$ coincides with
$\partial\Sigma$ then we calculate the integrand on level sets of
$f$ approaching $\partial\Sigma$ as follows.  Let $\Lambda$ be a
level set, $N$ the unit normal to $\Sigma$, $\nu$ the unit inner
normal to $\Lambda$ in the horizontal space (that is, in a
$t=\mathrm{const.}$ slice of $M\times\mathbb{R}$), and $\tau$ a
unit tangent to $\Sigma$ which is normal to $\Lambda$ (pointing
inside $\Sigma$), then a calculation [2] shows that
\begin{equation*}
H_{\Lambda,\overline{g}}=g_{\phi}(\tau,\nu)H_{\Lambda,g},
\end{equation*}
\begin{equation*}
\overline{g}(q,n_{\overline{g}})|_{\Lambda}=g_{\phi}(\tau,\nu)^{-1}
(H_{\Lambda,g}-g_{\phi}(N,\nu)\mathrm{Tr}_{\Lambda,g}k)
-H_{\Lambda,\overline{g}},
\end{equation*}
where $g_{\phi}=g+\phi^{2}dt^{2}$ is the metric on
$M\times\mathbb{R}$ and $H_{\Lambda,g}$,
$H_{\Lambda,\overline{g}}$ are the mean curvatures of $\Lambda$
with respect to $g$ and $\overline{g}$.  Thus since the Jang
surface $\Sigma$ blows-up to $\pm\infty$ at horizons
\begin{equation*}
g_{\phi}(\tau,\nu)\rightarrow 0,\text{ }\text{ }\text{ }\text{ }
g_{\phi}(N,\nu)\rightarrow\pm 1,\text{ }\text{ }\text{ as }\text{
}\text{ }\Lambda\rightarrow\partial\Sigma.
\end{equation*}
The apparent horizon equations
\begin{equation*}
H_{\partial\Sigma,g}\pm\mathrm{Tr}_{\partial\Sigma,g}k=0
\end{equation*}
then imply that the boundary integral at $S_{0}$ vanishes.  From
(5.3) we now have
\begin{equation}
M_{\mathrm{ADM}}\geq\sqrt{\frac{A_{\overline{g}}(S_{0})}{16\pi}}\geq
\sqrt{\frac{A_{g}(S_{0})}{16\pi}}\geq\sqrt{\frac{A}{16\pi}},
\end{equation}
after observing that $\overline{g}$ measures areas to be at least
as large as does $g$.\par
  Lastly we treat the case of equality.  By appealing to the above
arguments including (5.1), and using the full expression for
$\overline{R}$ in (2.3), we find that
\begin{equation*}
0\geq\int_{\widetilde{\Sigma}}\phi
\left(16\pi(\mu-|J|_{g})+|h-K|_{\Sigma}|_{\overline{g}}^{2}
+2|q|_{\overline{g}}^{2}\right)d\omega_{\overline{g}}.
\end{equation*}
As $\phi>0$ away from $\partial\widetilde{\Sigma}$ this implies
that
\begin{equation*}
\mu-|J|_{g}\equiv 0,\text{ }\text{ }\text{ }\text{
}h-K|_{\Sigma}\equiv 0,\text{ }\text{ }\text{ }\text{ }q\equiv 0,
\end{equation*}
from which we obtain $\overline{R}\equiv 0$.  Therefore
$(\widetilde{\Sigma},\overline{g})$ is isometric to
$(\mathbb{SC}^{3},g_{\mathrm{SC}})$ the exterior region of the
$t=0$ slice of the Schwarzchild spacetime, by the time symmetric
version of the PI.  Hence we may write
\begin{equation*}
\overline{g}=\left(1-\frac{2M_{\mathrm{ADM}}}{r}\right)^{-1}dr^{2}+r^{2}d\Omega^{2},
\end{equation*}
so that
\begin{equation*}
\phi=\sqrt{A_{\overline{g}}(S_{r})}H_{S_{r},\overline{g}}
=\sqrt{4\pi r^{2}}\sqrt{\overline{g}^{11}}\frac{2}{r}=4\sqrt{\pi}
\left(1-\frac{2M_{\mathrm{ADM}}}{r}\right)^{1/2}.
\end{equation*}
Moreover as $g=\overline{g}-\phi^{2}df^{2}$, it follows that the
graph map $G:M\rightarrow\mathbb{SC}^{4}$ given by $G(x)=(x,f(x))$
provides an isometric embedding of $(M-\widetilde{U},g)$ into the
Schwarzchild spacetime
$(\mathbb{SC}^{4},g_{\mathrm{SC}}-\phi^{2}dt^{2})$, where
$\partial \widetilde{U}$ is the image of $S_{0}$ in $M$. By (5.4)
(with all inequalities replaced by equalities) $A_{g}(\partial
\widetilde{U})=A$.  Then if $\partial U$ is the outermost minimal
area enclosure, we must have $\widetilde{U}\subset U$. Lastly a
calculation (see [2]) shows that $h-K|_{\Sigma}\equiv 0$ implies
that the extrinsic curvature of
$G(M-\widetilde{U})\subset\mathbb{SC}^{4}$ is given by $k$.\par
  Now some comments concerning the above methods.  The
definition of $\phi$ in (5.2) is not as simple as it appears, in
that $\phi$ is also present on the right-hand side due to the
definition of $\overline{g}$.  Furthermore the IMCF in $\Sigma$
depends on $f$ and $\phi$ for the same reason.  Thus we are not
only concerned with solving a single equation, namely the
generalized Jang equation (2.1), rather we must solve a system of
equations. We may write this system down in the following way. If
we take the level set formulation of IMCF (as in [6]), so that the
flow surfaces $S_{r}$ are given by the level sets $r=u(x)$ for
some function $u$ on $M$, then $u$ must satisfy the equation
\begin{equation}
\mathrm{div}_{\overline{g}}\left(\frac{\nabla_{\overline{g}}u}{|\nabla_{\overline{g}}u|}
\right)=|\nabla_{\overline{g}}u|
\end{equation}
in which the left-hand side represents the mean curvature of
$S_{r}$ and the right-hand side is the inverse speed of the flow.
It then follows that ([6])
\begin{equation*}
A_{\overline{g}}(S_{r})=e^{r}A_{\overline{g}}(S_{0})=e^{u}A_{\overline{g}}(S_{0}),
\end{equation*}
so by definition of $\phi$ and (5.5) we have
\begin{equation}
\phi^{2}=A_{\overline{g}}(S_{0})e^{u}|\nabla_{\overline{g}}u|^{2}.
\end{equation}
Equations (2.1), (5.5), and (5.6) now form a $3\times 3$
degenerate elliptic system for the unknowns $u$, $\phi$, and $f$.
In section $\S 3$ we have successfully solved this system for the
special case of spherically symmetric initial data, and have found
that the solution has the same behavior as conjectured in this
paper for the general case.  As each of the equations (2.1) and
(5.5) already have full existence theories in the classical case
when $\phi\equiv 1$ ([6], [13]), it is possible that similar
techniques will yield the corresponding theory for this
generalized Jang/IMCF system.\par
   Finally we mention that the method proposed here, which in a
nut shell can be thought of as simply integrating away the ``bad"
term from the expression of $\overline{R}$ in (2.3), can also be
modified to generalize the other known proof of the time symmetric
PI, namely the conformal flow proof of Bray [1].  When this is
done a new modified Jang/conformal flow system is generated. While
this method of proof would yield a stronger result (since it
applies to multiple black holes), the system obtained is less
tractable at the moment [2].

\begin{center}
\textbf{6.  Appendix A}
\end{center} \setcounter{equation}{0}
\setcounter{section}{6}

   In this appendix we confirm Theorem 1.  The following notation
will be used.  Suppose that $\Sigma$ is a smooth hypersurface
inside the warped product space
$(M\times\mathbb{R},g_{\phi}=g+\phi^{2}dt^{2})$, and let $e_{1}$,
$e_{2}$, $e_{3}$, $e_{4}$ be a local orthonormal frame for
$\Sigma$ with $e_{4}=N$ normal and $e_{1}$, $e_{2}$, $e_{3}$
tangent to $\Sigma$.  The Levi-Civita connection for $g_{\phi}$
will be denoted by $\nabla_{a}=\nabla_{e_{a}}$ and that for
$\overline{g}$, the induced metric on $\Sigma$, by
$\overline{\nabla}_{a}=\overline{\nabla}_{e_{a}}$.  Furthermore
$x^{1}$, $x^{2}$, $x^{3}$ will be local coordinates on $M$ with
$x^{4}=t$, and the second fundamental form of $\Sigma$ will be
denoted
\begin{equation*}
h_{ij}=h(e_{i},e_{j})=\langle\nabla_{e_{j}}N,e_{i}\rangle,
\end{equation*}
where $\langle\cdot,\cdot\rangle$ is the inner product of
$g_{\phi}$.\par
   The preliminary calculations will generalize those of [13].  We
have
\begin{equation*}
\overline{\nabla}_{j}\langle\partial_{x^{4}},N\rangle=\langle\nabla_{j}\partial_{x^{4}},N\rangle
+\langle\partial_{x^{4}},e_{i}\rangle h_{ij},
\end{equation*}
where the repeated index $i$ is summed from 1 to 3.  Next
\begin{eqnarray*}
\overline{\nabla}_{l}\overline{\nabla}_{j}\langle\partial_{x^{4}},N\rangle&=&
\langle\nabla_{l}\nabla_{j}
\partial_{x^{4}},N\rangle+\langle\nabla_{j}\partial_{x^{4}},e_{i}\rangle h_{il}
+\langle\partial_{x^{4}},e_{i}\rangle\overline{\nabla}_{l}h_{ij}\\
& &+\langle\nabla_{l}\partial_{x^{4}}, e_{i}\rangle
h_{ij}+\langle\partial_{x^{4}},\nabla_{l}e_{i}\rangle h_{ij}
-\langle\partial_{x^{4}},e_{p}\rangle\overline{\Gamma}^{p}_{li}h_{ij},
\end{eqnarray*}
where $\overline{\Gamma}_{li}^{p}$ are Christoffel symbols for
$\overline{g}$.  However the Codazzi equations give
\begin{equation*}
\overline{\nabla}_{l}h_{ij}-\overline{\nabla}_{i}h_{lj}=R_{Njil}
\end{equation*}
where $R_{Njil}$ are components of the Riemann tensor for
$g_{\phi}$, and we also have
\begin{eqnarray*}
\langle\partial_{x^{4}},\nabla_{l}e_{i}\rangle&=&
\langle\langle\partial_{x^{4}},N\rangle N
+\langle\partial_{x^{4}},e_{p}\rangle e_{p},\nabla_{l}e_{i}\rangle\\
&=&-\langle\partial_{x^{4}},N\rangle
h_{il}+\langle\partial_{x^{4}},e_{p}\rangle\overline{\Gamma}_{li}^{p}.
\end{eqnarray*}
Therefore, adopting the convention that indices $i$, $j$, $l$, $p$
run from 1 to 3 and $a$, $b$ (appearing later) run from 1 to 4, it
follows that
\begin{eqnarray*}
\Delta_{\overline{g}}\langle\partial_{x^{4}},N\rangle&=&\sum_{i}\overline{\nabla}_{i}
\overline{\nabla}_{i}\langle\partial_{x^{4}},N\rangle\\
&=& \sum_{i}\langle\nabla_{i}\nabla_{i}\partial_{x^{4}},N\rangle+
2\sum_{i,j}\langle\nabla_{j}\partial_{x^{4}},e_{i}\rangle h_{ij}\\
&
&+\sum_{i}\overline{\nabla}_{i}H\langle\partial_{x^{4}},e_{i}\rangle
+\sum_{i,j}R_{Njij}\langle\partial_{x^{4}},e_{i}\rangle
-|h|^{2}\langle\partial_{x^{4}},N\rangle,
\end{eqnarray*}
where
\begin{equation*}
H=\sum_{i}h_{ii},\text{ }\text{ }\text{ }\text{ }\text{
}|h|^{2}=\sum_{i,j} h_{ij}h_{ij}.
\end{equation*}
Moreover
\begin{equation*}
\overline{\nabla}_{i}H\langle\partial_{x^{4}},e_{i}\rangle=\nabla_{\partial_{x^{4}}}H
-\nabla_{N}H\langle\partial_{x^{4}},N\rangle=
-N(H)\langle\partial_{x^{4}},N\rangle,
\end{equation*}
and
\begin{equation*}
R_{Njil}\langle\partial_{x^{4}},e_{i}\rangle=-R_{NjNl}\langle\partial_{x^{4}},N\rangle
+\mathrm{Riem}(N,e_{j},\partial_{x^{4}},e_{l}).
\end{equation*}
Hence
\begin{eqnarray}
\Delta_{\overline{g}}\langle\partial_{x^{4}},N\rangle\!\!\!&=&\!\!\!
-(|h|^{2}+N(H)+\sum_{i}R_{NiNi})\langle\partial_{x^{4}},N\rangle
+\sum_{i}\langle\nabla_{i}\nabla_{i}\partial_{x^{4}},N\rangle\nonumber\\
& &\!\!\!+\sum_{i}\mathrm{Riem}(N,e_{i},\partial_{x^{4}},e_{i})
+2\sum_{i,j}\langle\nabla_{j}\partial_{x^{4}},e_{i}\rangle h_{ij}.
\end{eqnarray}\par
   Let $K$ be the extended version (by (2.2)) of the initial data.
We then define the extended versions of the local energy and
current densities by
\begin{equation*}
2\mu^{\mathrm{ext}}=R_{g_{\phi}}-\sum_{a,b}K_{ab}^{2}+(\sum_{a}K_{aa})^{2},\text{
}\text{ }\text{ }\text{ }J^{\mathrm{ext}}(e_{b})=\sum_{a}
(\nabla_{a}K_{ab}-\nabla_{b}K_{aa}),
\end{equation*}
where $R_{g_{\phi}}$ is the scalar curvature of $g_{\phi}$ and
$K_{ab}=K(e_{a},e_{b})$.  Notice that
\begin{equation*}
R_{g_{\phi}}=2\sum_{i}R_{NiNi}+\sum_{i,j}R_{ijij}
\end{equation*}
and by the Gauss equations
\begin{equation*}
\overline{R}_{ijpl}=R_{ijpl}+h_{ip}h_{jl}-h_{il}h_{jp},
\end{equation*}
(here $\overline{R}_{ijpl}$ denotes the Riemann tensor of
$\overline{g}$) which implies
\begin{equation*}
R_{g_{\phi}}=2\sum_{i}R_{NiNi}+\overline{R}-H^{2}+|h|^{2},
\end{equation*}
where $\overline{R}$ is the scalar curvature of $\overline{g}$. So
by definition of $\mu^{\mathrm{ext}}$,
\begin{equation*}
\sum_{i}R_{NiNi}=\mu^{\mathrm{ext}}+\frac{1}{2}\left(-\overline{R}+\sum_{a,b}K_{ab}^{2}
-(\sum_{a}K_{aa})^{2}-|h|^{2}+H^{2}\right).
\end{equation*}
Thus (6.1) becomes
\begin{eqnarray}
&
&\!\!\Delta_{\overline{g}}\langle\partial_{x^{4}},N\rangle\\
\!\!&=&\!\!-\left(\mu^{\mathrm{ext}}
+\frac{1}{2}|h|^{2}+\frac{1}{2}H^{2}+N(H)-\frac{1}{2}\overline{R}
+\frac{1}{2}\sum_{a,b}K_{ab}^{2}-\frac{1}{2}(\sum_{a}K_{aa})^{2}\right)
\langle\partial_{x^{4}},N\rangle\nonumber\\
&
&\!\!+\sum_{i}\mathrm{Riem}(N,e_{i},\partial_{x^{4}},e_{i})+\sum_{i}
\langle\nabla_{i}\nabla_{i}\partial_{x^{4}},N\rangle+2\sum_{i,j}\langle\nabla_{j}\partial_{x^{4}},
e_{i}\rangle h_{ij}.\nonumber
\end{eqnarray}\par
   We now obtain another expression for
$\Delta_{\overline{g}}\langle\partial_{x^{4}},N\rangle$.  First
extend the second fundamental form tensor $h$ to all of
$M\times\mathbb{R}$ by
\begin{equation*}
h(X,Y)=\langle\nabla_{Y}N,X\rangle,\text{ }\text{ }\text{ }\text{
} X,Y\in T_{x_{0}}(M\times\mathbb{R}),
\end{equation*}
so that
\begin{equation*}
h_{iN}=h(e_{i},N),\text{ }\text{ }\text{ }\text{
}h_{Ni}=h_{NN}=0,\text{ }\text{ }\text{ }\text{ }i=1,2,3.
\end{equation*}
Observe that
\begin{equation*}
\nabla_{\partial_{x^{4}}}N=\sum_{i,j}\langle\partial_{x^{4}},e_{i}\rangle
h_{ij}e_{j} +\sum_{j}\langle\partial_{x^{4}},N\rangle h_{jN}e_{j}
\end{equation*}
so that
\begin{equation}
h_{jN}=-\overline{\nabla}_{j}\log\langle\partial_{x^{4}},N\rangle
+\langle\partial_{x^{4}},N\rangle^{-1}(\langle\nabla_{j}\partial_{x^{4}},N\rangle
+\langle\nabla_{\partial_{x^{4}}}N, e_{j}\rangle),
\end{equation}
which implies
\begin{eqnarray*}
& &\!\!\Delta_{\overline{g}}\log\langle\partial_{x^{4}},N\rangle\\
&=&\!\!
-\langle\partial_{x^{4}},N\rangle^{-2}\sum_{j}(\overline{\nabla}_{j}\langle\partial_{x^{4}},N\rangle)^{2}
+\langle\partial_{x^{4}},N\rangle^{-1}\Delta_{\overline{g}}\langle\partial_{x^{4}},N\rangle\\
&=&\!\!-\sum_{j}[h_{jN}-\langle\partial_{x^{4}},N\rangle^{-1}(\langle\nabla_{j}\partial_{x^{4}},N\rangle
+\langle\nabla_{\partial_{x^{4}}}N,e_{j}\rangle)]^{2}\\
&
&\!\!+\langle\partial_{x^{4}},N\rangle^{-1}\Delta_{\overline{g}}\langle\partial_{x^{4}},N\rangle.
\end{eqnarray*}
With the help of (6.3) we have
\begin{eqnarray*}
&
&\!\!\langle\partial_{x^{4}},N\rangle^{-1}\Delta_{\overline{g}}\langle\partial_{x^{4}},N\rangle\\
&=&\!\!\sum_{j}[h_{jN}-\langle\partial_{x^{4}},N\rangle^{-1}(\langle\nabla_{j}\partial_{x^{4}},N\rangle
+\langle\nabla_{\partial_{x^{4}}}N,e_{j}\rangle)]^{2}\\
& &\!\!-\sum_{j}\overline{\nabla}_{j}h_{jN}
+\sum_{j}\overline{\nabla}_{j}\left(\frac{\langle\nabla_{j}\partial_{x^{4}},N\rangle}{\langle\partial_{x^{4}},N\rangle}
+\frac{\langle\nabla_{\partial_{x^{4}}}N,e_{i}\rangle}{\langle\partial_{x^{4}},N\rangle}\right)\\
&=&\!\!\sum_{j}(h_{jN}^{2}-2\langle\partial_{x^{4}},N\rangle^{-1}
\langle\nabla_{j}\partial_{x^{4}},N\rangle h_{jN}
+\langle\partial_{x^{4}},N\rangle^{-2}\langle\nabla_{\partial_{x^{4}}}N,e_{j}\rangle^{2})\\
&
&\!\!-\sum_{j}2\langle\partial_{x^{4}},N\rangle^{-1}\langle\nabla_{\partial_{x^{4}}}N,e_{j}\rangle
(h_{jN}-\langle\partial_{x^{4}},N\rangle^{-1}\langle\nabla_{j}\partial_{x^{4}},N\rangle)\\
&
&\!\!+\sum_{j}[\langle\partial_{x^{4}},N\rangle^{-1}(\langle\nabla_{j}\nabla_{j}\partial_{x^{4}},N\rangle
+\sum_{i}\langle\nabla_{j}\partial_{x^{4}},e_{i}\rangle h_{ji})-\overline{\nabla}_{j}h_{jN}]\\
&
&\!\!+\sum_{j}\left(\overline{\nabla}_{j}\frac{\langle\nabla_{\partial_{x^{4}}}N,e_{j}\rangle}
{\langle\partial_{x^{4}},N\rangle}
-\sum_{i}\langle\partial_{x^{4}},N\rangle^{-2}\langle\partial_{x^{4}},e_{i}\rangle
h_{ji}\langle\nabla_{j}\partial_{x^{4}},N\rangle \right).
\end{eqnarray*}
However
\begin{equation*}
\sum_{i}\langle\partial_{x^{4}},e_{i}\rangle
h_{ji}=-\langle\partial_{x^{4}},N\rangle h_{jN}
+\langle\nabla_{\partial_{x^{4}}}N,e_{j}\rangle
\end{equation*}
and
\begin{eqnarray*}
&
&\!\!\sum_{j}\overline{\nabla}_{j}\frac{\langle\partial_{x^{4}}N,e_{j}\rangle}
{\langle\partial_{x^{4}},N\rangle}\\
&=&\!\!\sum_{j}\langle\partial_{x^{4}},N\rangle^{-1}\overline{\nabla}_{j}
\langle\nabla_{\partial_{x^{4}}}N,e_{j}\rangle\\
&
&\!\!-\sum_{j}\langle\partial_{x^{4}},N\rangle^{-2}\langle\nabla_{\partial_{x^{4}}}N,e_{j}\rangle
(\langle\nabla_{j}\partial_{x^{4}},N\rangle-\langle\partial_{x^{4}},N\rangle
h_{jN} +\langle\nabla_{\partial_{x^{4}}}N,e_{j}\rangle),
\end{eqnarray*}
therefore
\begin{eqnarray}
&
&\!\!\langle\partial_{x^{4}},N\rangle^{-1}\Delta_{\overline{g}}\langle\partial_{x^{4}},N\rangle\\
&=&\!\!\sum_{j}(h_{jN}^{2}-\langle\partial_{x^{4}},N\rangle^{-1}
\langle\nabla_{j}\partial_{x^{4}},N\rangle h_{jN}
-\overline{\nabla}_{j}h_{jN})\nonumber\\
&
&\!\!+\sum_{j}\langle\partial_{x^{4}},N\rangle^{-1}(\langle\nabla_{j}\nabla_{j}\partial_{x^{4}},N\rangle
+\sum_{i}\langle\nabla_{j}\partial_{x^{4}},e_{i}\rangle h_{ji})\nonumber\\
&
&\!\!+\langle\partial_{x^{4}},N\rangle^{-1}\sum_{j}(\overline{\nabla}_{j}
\langle\nabla_{\partial_{x^{4}}}N,e_{j}\rangle-\langle\nabla_{\partial_{x^{4}}}N,e_{j}\rangle
h_{jN}).\nonumber
\end{eqnarray}\par
   In order to compare quantities appearing in (6.2) and (6.4) to the local
current density, we employ a formula on page 239 of [13]:
\begin{equation}
J^{\mathrm{ext}}(N)=\sum_{i}\overline{\nabla}_{i}K_{i4}-N(\sum_{i}K_{ii})
+K_{NN}H-\sum_{i,j}K_{ij}h_{ij}-2\sum_{i}K_{iN}h_{iN}.
\end{equation}
This formula still remains valid in our situation.  To see this
observe that
\begin{equation*}
J^{\mathrm{ext}}(N)=J^{\mathrm{ext}}(e_{4})=\sum_{i}(\nabla_{i}K_{iN}-\nabla_{N}K_{ii}).
\end{equation*}
Moreover if
$(\overline{\delta}^{ij})=(\overline{g}(e_{i},e_{j}))^{-1}$ then
\begin{eqnarray*}
\sum_{i}\nabla_{N}K_{ii}&=&\overline{\delta}^{ij}(N(K_{ij})-2\Gamma_{Ni}^{a}K_{ja})\\
&=&N(\sum_{i}K_{ii})-N(\overline{\delta}^{ij})K_{ij}-2\sum_{i,j}\Gamma_{iN}^{j}K_{ij}
-2\sum_{i}\Gamma_{Ni}^{N}K_{iN}\\
&=&N(\sum_{i}K_{ii})+\sum_{i}2h_{iN}K_{iN}
\end{eqnarray*}
since
\begin{equation*}
\Gamma_{Ni}^{N}=\langle
N,\nabla_{N}e_{i}\rangle=-\langle\nabla_{N}N,e_{i}\rangle=-h_{iN},
\end{equation*}
\begin{equation*}
\Gamma_{iN}^{j}=\langle e_{j},\nabla_{i}N\rangle=h_{ij},\text{
}\text{ }\text{ }\text{ }N(\overline{\delta}^{ij})=-2h_{ij},
\end{equation*}
and
\begin{eqnarray*}
\sum_{i}\nabla_{i}K_{iN}&=&
\sum_{i}(e_{i}(K_{iN})-K(\nabla_{i}e_{i},N)-K(e_{i},\nabla_{i}N))\\
&=&\sum_{i}\overline{\nabla}_{i}K_{iN}+HK_{NN}-\sum_{i,j}h_{ij}K_{ij}
\end{eqnarray*}
since
\begin{equation*}
\nabla_{i}e_{i}=\Gamma_{ii}^{N}N+\sum_{j}\Gamma_{ii}^{j}e_{j}
=-h_{ii}N+\sum_{j}\overline{\Gamma}_{ii}^{j}e_{j},\text{ }\text{
}\text{ }\text{ }\nabla_{i}N=\sum_{j}h_{ij}e_{j}.
\end{equation*}
The desired formula now follows.\par
  Equations (6.2) and (6.4) yield an expression for $\mu^{\mathrm{ext}}$.
Then by combining this expression with (6.5) we arrive at
\begin{eqnarray}
& &2(\mu^{\mathrm{ext}}\!-\!J^{\mathrm{ext}}(N))\\
&=&\!\!\!
\overline{R}-\sum_{i,j}(h_{ij}-K_{ij})^{2}-2\sum_{i}(h_{iN}-K_{iN})^{2}
+2\sum_{i}\overline{\nabla}_{i}(h_{iN}-K_{iN})\nonumber\\
&
&\!\!\!+(\sum_{i}K_{ii})^{2}\!-\!H^{2}+2K_{NN}(\sum_{i}K_{ii}\!-\!H)
+2N(\sum_{i}K_{ii}\!-\!H)\nonumber\\
&
&\!\!\!+2\langle\partial_{x^{4}},N\rangle^{-1}\sum_{i}(\mathrm{Riem}(N,e_{i},\partial_{x^{4}},
e_{i})+\langle\nabla_{i}\partial_{x^{4}},\sum_{j}h_{ij}e_{j}+h_{iN}N\rangle\nonumber\\
& &\!\!\!+\langle\nabla_{\partial_{x^{4}}}N,e_{i}\rangle h_{iN}
-\overline{\nabla}_{i}\langle\nabla_{\partial_{x^{4}}}N,e_{i}\rangle),\nonumber
\end{eqnarray}
where the repeated indices $i$, $j$ are summed from 1 to 3.\par
   The remainder of the proof will consist of evaluating certain terms
from (6.6) in local coordinates.  We assume from now on that
$\Sigma$ satisfies the generalized Jang equation (2.1), so that
the 5th, 6th, and 7th terms on the right-hand side of (6.6)
vanish. We also assume that $\Sigma$ is given by the graph of a
function $x^{4}=f(x^{1},x^{2},x^{3})$, and we will write
$f_{,i}=\partial f/\partial x^{i}$, $f^{i}=g^{ij}f_{,j}$.  Let
\begin{equation*}
X_{i}=\partial_{x^{i}}+f_{,i}\partial_{x^{4}},\text{ }\text{
}\text{ }\text{ }i=1,2,3,
\end{equation*}
be tangent vectors to $\Sigma$ and
\begin{equation*}
N=\frac{f^{i}\partial_{x^{i}}-\phi^{-2}\partial_{x^{4}}}
{\sqrt{\phi^{-2}+|\nabla_{g}f|^{2}}}
\end{equation*}
be the unit normal to $\Sigma$.  Also we will write
\begin{equation*}
\overline{g}_{ij}=\overline{g}(X_{i},X_{j})=g_{ij}+\phi^{2}
f_{,i}f_{,j},\text{ }\text{ }\text{ }\text{
}\overline{g}^{ij}=g^{ij}-\frac{f^{i}f^{j}}{\phi^{-2}+|\nabla_{g}f|^{2}}.
\end{equation*}
The next three claims will simplify (6.6).\medskip

\textbf{Claim 1.}
\begin{equation*}
\langle\partial_{x^{4}},N\rangle^{-1}\sum_{i}\mathrm{Riem}(N,e_{i},\partial_{x^{4}},e_{i})
=-\phi^{-1}\Delta_{g}\phi
\end{equation*}

\textit{Proof.} Christoffel symbols for the metric $g_{\phi}$ in
the above local coordinates are given by
\begin{equation}
\widehat{\Gamma}_{44}^{4}=\widehat{\Gamma}_{ij}^{4}=\widehat{\Gamma}_{i4}^{j}=0,\text{
}\text{ }\text{ }\text{ }1\leq i,j\leq 3,
\end{equation}
\begin{equation*}
\widehat{\Gamma}_{i4}^{4}=(\log\phi)_{,i},\text{ }\text{ }\text{
}\text{ }\widehat{\Gamma}_{44}^{i}=-\phi\phi^{i}.
\end{equation*}
Note that $\widehat{\Gamma}_{ij}^{k}$ are the Christoffel symbols
for the initial data metric $g$ when $1\leq i,j,k\leq 3$.  The
Riemann tensor is then given by
\begin{eqnarray}
\widehat{R}_{4ijk}&=&(g_{\phi})_{4a}\widehat{R}^{a}_{ijk}=\phi^{2}\widehat{R}^{4}_{ijk}\\
&=&\phi^{2}(\widehat{\Gamma}_{ik,j}^{4}-\widehat{\Gamma}_{ij,k}^{4}
+\widehat{\Gamma}_{ik}^{b}\widehat{\Gamma}_{bj}^{4}
-\widehat{\Gamma}_{ij}^{b}\widehat{\Gamma}_{bk}^{4})\nonumber\\
&=&0,\nonumber
\end{eqnarray}
\begin{eqnarray}
\widehat{R}_{4i4j}&=&(g_{\phi})_{4a}\widehat{R}^{a}_{i4j}=\phi^{2}\widehat{R}^{4}_{i4j}\\
&=&\phi^{2}(\widehat{\Gamma}_{ij,4}^{4}-\widehat{\Gamma}_{i4,j}^{4}
+\widehat{\Gamma}_{ij}^{b}\widehat{\Gamma}_{b4}^{4}
-\widehat{\Gamma}_{i4}^{b}\widehat{\Gamma}_{bj}^{4})\nonumber\\
&=&-\phi\phi_{;ij},\nonumber
\end{eqnarray}
where the semicolon denotes covariant differentiation with respect
to $g$.  Moreover if $1\leq i,j,k,l\leq 3$ then
$\widehat{R}_{ijkl}$ are just the components of the Riemann tensor
for $g$.  Therefore with the help of (6.8) and (6.9) we have
\begin{eqnarray*}
\sum_{i}\mathrm{Riem}(N,e_{i},\partial_{x^{4}},e_{i})
&=&\overline{g}^{ij}
\mathrm{Riem}\left(\frac{f^{l}\partial_{x^{l}}-\phi^{-2}\partial_{x^{4}}}
{\langle\partial_{x^{4}},N\rangle^{-1}},X_{i},\partial_{x^{4}},X_{j}\right)\\
&=&\langle\partial_{x^{4}},N\rangle\overline{g}^{ij}(\phi^{-2}\widehat{R}_{4i4j}
-f^{l}\widehat{R}_{li4j}-f^{l}f_{,i}\widehat{R}_{l44j})\\
&=&-\langle\partial_{x^{4}},N\rangle\phi^{-1}\Delta_{g}\phi.
\end{eqnarray*}
Q.E.D.\medskip

\textbf{Claim 2.}
\begin{eqnarray*}
&
&\langle\partial_{x^{4}},N\rangle^{-1}\sum_{i}\langle\nabla_{i}\partial_{x^{4}},\sum_{j}h_{ij}e_{j}+h_{iN}N\rangle\\
&
&=-|\nabla_{\overline{g}}\log\phi+\phi\phi^{l}f_{,l}\nabla_{\overline{g}}f|^{2}
-\langle\partial_{x^{4}},N\rangle^{-1}\overline{g}^{ij}f_{,j}\phi\phi^{l}h(X_{i},X_{l})
\end{eqnarray*}

\textit{Proof.}  First observe that
\begin{equation*}
\langle\nabla_{i}\partial_{x^{4}},\sum_{j}h_{ij}e_{j}+h_{iN}N\rangle
=h(e_{i},\nabla_{e_{i}}\partial_{x^{4}}),
\end{equation*}
and therefore
\begin{equation*}
\sum_{i}\langle\nabla_{i}\partial_{x^{4}},\sum_{j}h_{ij}e_{j}+h_{iN}N\rangle
=\overline{g}^{ij}h(X_{i},\nabla_{X_{j}}\partial_{x^{4}}).
\end{equation*}
Now compute with the help of (6.7):
\begin{eqnarray*}
h(X_{i},\nabla_{X_{j}}\partial_{x^{4}})&=&
h(X_{i},\widehat{\Gamma}_{j4}^{a}\partial_{x^{a}}
+f_{,j}\widehat{\Gamma}_{44}^{a}\partial_{x^{a}})\\
&=&((\log\phi)_{,j}+f_{,j}f_{,k}\phi\phi^{k})h(X_{i},\partial_{x^{4}})
-f_{,j}\phi\phi^{k}h(X_{i},X_{k}).
\end{eqnarray*}
Moreover
\begin{eqnarray}
h(X_{i},\partial_{x^{4}})&=&-\langle N,\nabla_{\partial_{x^{4}}}X_{i}\rangle\\
&=&-\langle
N,\widehat{\Gamma}_{4i}^{a}\partial_{x^{a}}+f_{,i}\widehat{\Gamma}_{44}^{a}
\partial_{x^{a}}\rangle\nonumber\\
&=&-\langle\partial_{x^{4}},N\rangle((\log\phi)_{,i}+f_{,i}f_{,k}\phi\phi^{k}).\nonumber
\end{eqnarray}
Q.E.D.\medskip

\textbf{Claim 3.}
\begin{eqnarray*}
&
&\!\!\!\!\langle\partial_{x^{4}},N\rangle^{-1}\sum_{i}(\langle\nabla_{\partial_{x^{4}}}N,e_{i}\rangle
h_{iN}-
\overline{\nabla}_{i}\langle\nabla_{\partial_{x^{4}}}N,e_{i}\rangle)\\
&=&\!\!\!\! \phi^{-1}\Delta_{g}\phi
+\langle\partial_{x^{4}},N\rangle^{-1}\overline{g}^{ij}f_{,j}\phi\phi^{l}
h(X_{i},X_{l})+|\nabla_{\overline{g}}\log\phi
+\phi\phi^{l}f_{,l}\nabla_{\overline{g}}f|^{2}
\end{eqnarray*}

\textit{Proof.}  First note that
\begin{eqnarray*}
\nabla_{\partial_{x^{4}}}N&=&\frac{1}{\sqrt{\phi^{-2}+|\nabla_{g}f|^{2}}}
(f^{i}\widehat{\Gamma}_{i4}^{a}\partial_{x^{a}}-\phi^{-2}\widehat{\Gamma}_{44}^{a}
\partial_{x^{a}})\\
&=&-\langle\partial_{x^{4}},N\rangle(\log\phi)^{i}X_{i},
\end{eqnarray*}
and therefore
\begin{eqnarray}
&
&\langle\partial_{x^{4}},N\rangle^{-1}\sum_{i}\langle\nabla_{\partial_{x^{4}}}N,e_{i}\rangle h_{iN}\\
&=&\langle\partial_{x^{4}},N\rangle^{-1}\overline{g}^{ij}\langle\nabla_{\partial_{x^{4}}}N,X_{j}\rangle
h(X_{i},N)\nonumber\\
&=&-\overline{g}^{ij}(\log\phi)^{k}\overline{g}_{kj}
h\left(X_{i},\frac{f^{l}\partial_{x^{l}}-\phi^{-2}\partial_{x^{4}}}
{\sqrt{\phi^{-2}+|\nabla_{g}f|^{2}}}\right)\nonumber\\
&=&-\frac{1}{\sqrt{\phi^{-2}+|\nabla_{g}f|^{2}}}(\log\phi)^{i}
f^{j}h(X_{i},X_{j})+\sqrt{\phi^{-2}+|\nabla_{g}f|^{2}}(\log\phi)^{i}
h(X_{i},\partial_{x^{4}}).\nonumber
\end{eqnarray}\par
  For the other term we have
\begin{eqnarray}
&
&-\langle\partial_{x^{4}},N\rangle^{-1}\sum_{i}\overline{\nabla}_{i}\langle\partial_{x^{4}}N,e_{i}\rangle\\
&=&-\langle\partial_{x^{4}},N\rangle^{-1}\overline{g}^{ij}
\overline{\nabla}_{X_{i}}\langle\partial_{x^{4}}N,X_{j}\rangle\nonumber\\
&=&\sqrt{\phi^{-2}+|\nabla_{g}f|^{2}}\overline{g}^{ij}\overline{\nabla}_{X_{i}}
\left(\frac{(\log\phi)^{k}\overline{g}_{kj}}{\sqrt{\phi^{-2}+|\nabla_{g}f|^{2}}}\right)\nonumber\\
&=&\overline{\nabla}_{X_{i}}(\log\phi)^{i}
+\frac{\phi^{-2}|\nabla_{g}\log\phi|^{2}}{\phi^{-2}+|\nabla_{g}f|^{2}}
-\frac{(\log\phi)^{i}f^{j}f_{;ij}}{\phi^{-2}+|\nabla_{g}f|^{2}}.\nonumber
\end{eqnarray}
Also
\begin{eqnarray}
\overline{\nabla}_{X_{i}}(\log\phi)^{i}&=&X_{i}[(\log\phi)^{i}]
+\overline{\Gamma}_{ik}^{i}(\log\phi)^{k}\\
&=&\partial_{x^{i}}[g^{il}(\log\phi)_{,l}]
+\widehat{\Gamma}_{ik}^{i}(\log\phi)^{k}
+(\overline{\Gamma}_{ik}^{i}-\widehat{\Gamma}_{ik}^{i})
(\log\phi)^{k}\nonumber\\
&=&\Delta_{g}\log\phi+(\overline{\Gamma}_{ik}^{i}-\widehat{\Gamma}_{ik}^{i})
(\log\phi)^{k}\nonumber
\end{eqnarray}
since
\begin{equation*}
\partial_{x^{i}}g^{il}=-g^{ij}\widehat{\Gamma}_{jk}^{k}-g^{jk}\widehat{\Gamma}_{jk}^{i},
\end{equation*}
where the overline indicates Christoffel symbols for the induced
metric $\overline{g}$.  In order to calculate the difference of
Christoffel symbols appearing above, notice that
\begin{eqnarray*}
\overline{\Gamma}_{jk}^{l}X_{l}=\overline{\nabla}_{X_{j}}X_{k}
&=&\nabla_{X_{j}}X_{k}+h(X_{k},X_{j})N\\
&=&\nabla_{\partial_{x^{j}}}(\partial_{x^{k}}+f_{,k}\partial_{x^{4}})
+f_{,j}\nabla_{\partial_{x^{4}}}(\partial_{x^{k}}+f_{k}\partial_{x^{4}})
+h(X_{k},X_{j})N\\
&=&\left(\widehat{\Gamma}_{jk}^{l}-\phi\phi^{l}f_{,j}f_{,k}
+\frac{f^{l}h(X_{k},X_{j})}{\sqrt{\phi^{-2}+|\nabla_{g}f|^{2}}}\right)
\partial_{x^{l}}\\
& &+\left(f_{,jk}+(\log\phi)_{,j}f_{,k}+(\log\phi)_{,k}f_{,j}
-\frac{\phi^{-2}h(X_{k},X_{j})}{\sqrt{\phi^{-2}+|\nabla_{g}f|^{2}}}\right)
\partial_{x^{4}}\\
&=&\left(\widehat{\Gamma}_{jk}^{l}-\phi\phi^{l}f_{,j}f_{,k}
+\frac{f^{l}h(X_{k},X_{j})}{\sqrt{\phi^{-2}+|\nabla_{g}f|^{2}}}\right)
X_{l},
\end{eqnarray*}
where we have used the formula
\begin{eqnarray}
h(X_{i},X_{j})&=&\langle\nabla_{X_{j}}N,X_{i}\rangle\\
&=&\frac{1}{\sqrt{\phi^{-2}+|\nabla_{g}f|^{2}}}
(f_{;ij}+(\log\phi)_{,i}f_{,j}+(\log\phi)_{,j}f_{,i}
+\phi\phi^{l}f_{,l}f_{,i}f_{,j})\nonumber
\end{eqnarray}
which is easily established from (6.7).  Hence
\begin{equation}
\overline{\Gamma}_{jk}^{l}=\widehat{\Gamma}_{jk}^{l}
-\phi\phi^{l}f_{,j}f_{,k}
+\frac{f^{l}h(X_{k},X_{j})}{\sqrt{\phi^{-2}+|\nabla_{g}f|^{2}}}.
\end{equation}
By combining (6.12)-(6.15) we arrive at
\begin{equation}
-\langle\partial_{x^{4}},N\rangle^{-1}\sum_{i}\overline{\nabla}_{i}\langle\partial_{x^{4}}N,e_{i}\rangle
=\phi^{-1}\Delta_{g}\phi.
\end{equation}\par
  Lastly, with (6.11) and (6.16) the desired result is obtained
after making a short calculation (using (6.10)) to show that
\begin{equation*}
\sqrt{\phi^{-2}+|\nabla_{g}f|^{2}}(\log\phi)^{i}h(X_{i},\partial_{x^{4}})
=|\nabla_{\overline{g}}\log\phi+\phi\phi^{l}f_{,l}\nabla_{\overline{g}}f|^{2},
\end{equation*}
and also
\begin{equation*}
\overline{g}^{ij}f_{,j}=\frac{\phi^{-2}f^{i}}{\phi^{-2}+|\nabla_{g}f|^{2}}.
\end{equation*}
Q.E.D.\medskip
  Claims 1, 2, and 3 show that the last four terms of (6.6)
cancel to yield
\begin{eqnarray}
& &2(\mu^{\mathrm{ext}}-J^{\mathrm{ext}}(N))\\
&=&\overline{R}-\sum_{i,j}(h_{ij}-K_{ij})^{2}-2\sum_{i}(h_{iN}-K_{iN})^{2}
+2\sum_{i}\overline{\nabla}_{i}(h_{iN}-K_{iN}).\nonumber
\end{eqnarray}\par
  We continue by writing the ``extended" energy and current
densities in terms of the original densities.\medskip

\textbf{Claim 4.}
\begin{equation*}
\mu^{\mathrm{ext}}-J^{\mathrm{ext}}(N)=8\pi(\mu-J(w))
-\phi^{-1}\Delta_{g}\phi+Q(k,\phi,f)
\end{equation*}
\textit{where}
\begin{equation*}
w=\frac{f^{i}\partial_{x^{i}}}{\sqrt{\phi^{-2}+|\nabla_{g}f|^{2}}},
\end{equation*}
\textit{and}
\begin{equation*}
Q(k,\phi,f)=\phi^{-2}(\mathrm{Tr}_{g}k)k_{44}
+\frac{f^{i}}{\sqrt{\phi^{-2}+|\nabla_{g}f|^{2}}}(
\phi^{-2}k_{44,i}-\phi^{-2}(\log\phi)_{,i}k_{44}
-(\log\phi)^{j}k_{ij})
\end{equation*}
\textit{with the extension term $k_{44}$ given by (2.2) and
$k_{ij}=k(\partial_{x^{i}},\partial_{x^{j}})$.}\medskip

\textit{Proof.}  We first treat the energy density
\begin{eqnarray*}
2\mu^{\mathrm{ext}}&=&R_{g_{\phi}}-g_{\phi}^{ac}g_{\phi}^{bd}
K_{ab}K_{cd}+(g_{\phi}^{ab}K_{ab})^{2}\\
&=&R_{g_{\phi}}-R+16\pi\mu+2\phi^{-2}(\mathrm{Tr}_{g}k)k_{44},
\end{eqnarray*}
where $R$ is the scalar curvature of $g$.  Moreover by (6.7)
\begin{eqnarray*}
R_{g_{\phi}}&=&g_{\phi}^{jl}\widehat{\Gamma}_{jl,k}^{k}-g_{\phi}^{ik}
\widehat{\Gamma}_{ij,k}^{j}+g_{\phi}^{jl}\widehat{\Gamma}_{ik}^{k}
\widehat{\Gamma}_{jl}^{i}-g_{\phi}^{ik}\widehat{\Gamma}_{kl}^{j}
\widehat{\Gamma}_{ij}^{l}\\
&=&R+g_{\phi}^{44}\widehat{\Gamma}_{44,k}^{k}-g_{\phi}^{ik}
\widehat{\Gamma}_{i4,k}^{4}+g_{\phi}^{jl}\widehat{\Gamma}_{jl}^{i}
\widehat{\Gamma}_{i4}^{4}+g_{\phi}^{jl}\widehat{\Gamma}_{jl}^{4}\widehat{\Gamma}_{4k}^{k}\\
& &
+g_{\phi}^{44}\widehat{\Gamma}_{44}^{i}\widehat{\Gamma}_{ik}^{k}
-g_{\phi}^{ik}\widehat{\Gamma}_{k4}^{4}\widehat{\Gamma}_{i4}^{4}
-g_{\phi}^{44}\widehat{\Gamma}_{4l}^{j}\widehat{\Gamma}_{4j}^{l}\\
&=&R-2\phi^{-1}\Delta_{g}\phi,
\end{eqnarray*}
and therefore
\begin{equation*}
\mu^{\mathrm{ext}}=8\pi\mu-\phi^{-1}\Delta_{g}\phi+\phi^{-2}(\mathrm{Tr}_{g}k)k_{44}.
\end{equation*}\par
  Now consider the current density.  If $1\leq i\leq 3$ then
\begin{eqnarray*}
J^{\mathrm{ext}}(\partial_{x^{i}})&=&g_{\phi}^{ab}K_{bi;a}-g_{\phi}^{ab}
K_{ab;i}\\
&=&g_{\phi}^{ab}(\partial_{x^{a}}K_{bi}-\widehat{\Gamma}_{ab}^{c}K_{ci}
-\widehat{\Gamma}_{ia}^{c}K_{cb})\\
&
&-g_{\phi}^{ab}(\partial_{x^{i}}K_{ab}-\widehat{\Gamma}_{bi}^{c}K_{ca}
-\widehat{\Gamma}_{ai}^{c}K_{cb})\\
&=&8\pi J(\partial_{x^{i}})-\phi^{-2}k_{44,i}
+\phi^{-2}(\log\phi)_{,i}k_{44} +(\log\phi)^{j}k_{ij}.
\end{eqnarray*}
In addition
\begin{eqnarray*}
J^{\mathrm{ext}}(\partial_{x^{4}})&=&g_{\phi}^{ab}K_{b4;a}-g_{\phi}^{ab}
K_{ab;4}\\
&=&g_{\phi}^{ab}(\partial_{x^{a}}K_{b4}-\widehat{\Gamma}_{ab}^{c}K_{c4}
-\widehat{\Gamma}_{4a}^{c}K_{cb})\\
&
&-g_{\phi}^{ab}(\partial_{x^{4}}K_{ab}-\widehat{\Gamma}_{b4}^{c}K_{ca}
-\widehat{\Gamma}_{a4}^{c}K_{cb})\\
&=&0,
\end{eqnarray*}
so that
\begin{eqnarray*}
J^{\mathrm{ext}}(N)&=&J^{\mathrm{ext}}\left(\frac{f^{i}\partial_{x^{i}}
-\phi^{-2}\partial_{x^{4}}}{\sqrt{\phi^{-2}+|\nabla_{g}f|^{2}}}\right)\\
&=&8\pi J(w)-\frac{f^{i}}{\sqrt{\phi^{-2}+|\nabla_{g}f|^{2}}}(
\phi^{-2}k_{44,i}-\phi^{-2}(\log\phi)_{,i}k_{44}
-(\log\phi)^{j}k_{ij}).
\end{eqnarray*}
Q.E.D.\medskip

  Our next goal will be to simplify the expression for
$Q(k,\phi,f)$.  To this end we will need the following
\medskip

\textbf{Claim 5.}
\begin{equation*}
\sum_{i}(h_{iN}-K_{iN})^{2}= |q|_{\overline{g}}^{2}
+|\nabla_{\overline{g}}\log\phi|^{2}-\overline{g}^{ij}
q_{i}(\log\phi)_{,j},
\end{equation*}
\textit{where}
\begin{equation*}
q_{i}=\frac{f^{j}}{\sqrt{\phi^{-2}+|\nabla_{g}f|^{2}}}(
h(X_{i},X_{j})-K(X_{i},X_{j})).
\end{equation*}

\textit{Proof.}  Employ (6.10) and (2.2) to find
\begin{eqnarray}
&
&h(X_{i},N)-K(X_{i},N)\\
&=&\frac{1}{\sqrt{\phi^{-2}+|\nabla_{g}f|^{2}}}
(h(X_{i},f^{j}\partial_{x^{j}}-\phi^{-2}\partial_{x^{4}})
-K(X_{i},f^{j}\partial_{x^{j}}-\phi^{-2}\partial_{x^{4}}))\nonumber\\
&=&\frac{f^{j}}{\sqrt{\phi^{-2}+|\nabla_{g}f|^{2}}}
(h(X_{i},X_{j})-K(X_{i},X_{j}))+\sqrt{\phi^{-2}+|\nabla_{g}f|^{2}}(
f_{,i}k_{44}-h(X_{i},\partial_{x^{4}}))\nonumber\\
&=&q_{i}-(\log\phi)_{,i}.\nonumber
\end{eqnarray}
Q.E.D.\medskip

  As an immediate corollary of (6.18) we also have\medskip

\textbf{Claim 6.}
\begin{equation*}
\sum_{i}\overline{\nabla}_{i}(h_{iN}-K_{iN})=\mathrm{div}_{\overline{g}}q-\Delta_{\overline{g}}\log\phi
\end{equation*}

  We now come to the simplification of $Q(k,\phi,f)$.\medskip

\textbf{Claim 7.}
\begin{eqnarray*}
& &\sum_{i}(h_{i4}-K_{i4})^{2}+\sum_{i,j}(h_{ij}-K_{ij})^{2}+2Q(k,\phi,f)\\
&=&|h-K|_{\Sigma}|_{\overline{g}}^{2}+|q|_{\overline{g}}^{2}
-2\Delta_{\overline{g}}\log\phi
-|\nabla_{\overline{g}}\log\phi|^{2}
+2\phi^{-1}\Delta_{g}\phi\nonumber
\end{eqnarray*}

\textit{Proof.}  Using (2.2), the following term of $Q(k,\phi,f)$
(from Claim 4) may be calculated by
\begin{equation}
\frac{f^{i}(\log\phi)^{j}k_{ij}}{\sqrt{\phi^{-2}+|\nabla_{g}f|^{2}}}
=\frac{\phi^{-2}(\phi^{l}f_{,l})^{2}}{\phi^{-2}+|\nabla_{g}f|^{2}}.
\end{equation}
To see this observe that
\begin{eqnarray*}
f^{j}k_{ij}&=&K(\partial_{x^{i}}+f_{,i}\partial_{x^{4}},f^{j}\partial_{x^{j}})\\
&=&\sqrt{\phi^{-2}+|\nabla_{g}f|^{2}}K(X_{i},N)+\phi^{-2}f_{,i}k_{44},
\end{eqnarray*}
and
\begin{eqnarray*}
2(\log\phi)^{i}K(X_{i},N)&=&-2\langle\partial_{x^{4}},N\rangle^{-1}K(\nabla_{\partial_{x^{4}}}N,N)\\
&=&-\langle\partial_{x^{4}},N\rangle^{-1}(\partial_{x^{4}}K(N,N)-(\nabla_{\partial_{x^{4}}}K)(N,N))\\
&=&0.
\end{eqnarray*}
Therefore by calculating the remaining terms of $Q(k,\phi,f)$ in a
straight forward way, we have
\begin{eqnarray*}
Q(k,\phi,f)&=&\frac{(\mathrm{Tr}_{g}k)\phi^{-1}\phi^{l}f_{,l}}{\sqrt{\phi^{-2}+|\nabla_{g}f|^{2}}}
-\frac{\phi^{-2}(\phi^{l}f_{,l})^{2}}{\phi^{-2}+|\nabla_{g}f|^{2}}
+\frac{\phi^{-1}f^{i}f^{j}\phi_{;ij}}{\phi^{-2}+|\nabla_{g}f|^{2}}\\
& &+\frac{(\log\phi)^{i}f^{j}f_{;ij}}{\phi^{-2}+|\nabla_{g}f|^{2}}
+\frac{\phi^{-4}(\phi^{l}f_{,l})^{2}}{(\phi^{-2}+|\nabla_{g}f|^{2})^{2}}
-\frac{\phi^{-1}(\phi^{l}f_{,l})f^{i}f^{j}f_{;ij}}{(\phi^{-2}+|\nabla_{g}f|^{2})^{2}}.
\end{eqnarray*}
Moreover with the help of (6.15)
\begin{eqnarray*}
\frac{\phi^{-1}f^{i}f^{j}\phi_{;ij}}{\phi^{-2}+|\nabla_{g}f|^{2}}
&=&\phi^{-1}\Delta_{g}\phi-\phi^{-1}\overline{g}^{ij}\phi_{;ij}\\
&=&\phi^{-1}(\Delta_{g}\phi-\Delta_{\overline{g}}\phi)
+\phi^{-1}(\widehat{\Gamma}_{ij}^{k}-\overline{\Gamma}_{ij}^{k})\phi_{,k}\\
&=&\phi^{-1}(\Delta_{g}\phi-\Delta_{\overline{g}}\phi)
+|\nabla_{g}\phi|^{2}\overline{g}^{ij}f_{,i}f_{,j}
-\frac{\phi^{-1}(\phi^{l}f_{,l})H}{\sqrt{\phi^{-2}+|\nabla_{g}f|^{2}}},
\end{eqnarray*}
and
\begin{eqnarray*}
\mathrm{Tr}_{g}k&=&g^{ij}K(X_{i},X_{j})-|\nabla_{g}f|^{2}k_{44}\\
&=&\overline{g}^{ij}K(X_{i},X_{j})-|\nabla_{g}f|^{2}k_{44}\\
& &+\frac{f^{i}f^{j}}{\phi^{-2}+|\nabla_{g}f|^{2}}
(K(X_{i},X_{j})-h(X_{i},X_{j}))
+\frac{f^{i}f^{j}}{\phi^{-2}+|\nabla_{g}f|^{2}}h(X_{i},X_{j}).
\end{eqnarray*}
It follows that with (6.14)
\begin{eqnarray}
Q(k,\phi,f)&=&-\frac{\phi^{-1}(\phi^{l}f_{,l})}{\sqrt{\phi^{-2}+|\nabla_{g}f|^{2}}}
(H-\sum_{i}K_{ii})+\Delta_{g}\log\phi-\Delta_{\overline{g}}\log\phi\\
&
&-\frac{\phi^{-1}(\phi^{l}f_{,l})f^{i}f^{j}}{(\phi^{-2}+|\nabla_{g}f|^{2})^{3/2}}
(h(X_{i},X_{j})-K(X_{i},X_{j}))
+\frac{\phi^{-2}(\phi^{l}f_{,l})^{2}}{\phi^{-2}+|\nabla_{g}f|^{2}}\nonumber\\
& &+\frac{(\log\phi)^{i}f^{j}f_{;ij}}{\phi^{-2}+|\nabla_{g}f|^{2}}
+\frac{|\nabla_{g}f|^{2}|\nabla_{g}\log\phi|^{2}}{\phi^{-2}+|\nabla_{g}f|^{2}}.\nonumber
\end{eqnarray}\par
  Now set
\begin{equation*}
p_{ij}=h(X_{i},X_{j})-K(X_{i},X_{j}).
\end{equation*}
Many of the terms appearing in (6.20) are similar to those
appearing in the following expression, which is derived from
(6.18)
\begin{eqnarray}
& &\sum_{i,j}(h_{ij}-K_{ij})^{2}+\sum_{i}(h_{iN}-K_{iN})^{2}\\
&=&g^{il}g^{jk}p_{ij}p_{lk}-\frac{g^{il}f^{j}f^{k}}{\phi^{-2}+|\nabla_{g}f|^{2}}
p_{ij}p_{lk}-\frac{2f^{l}(\log\phi)^{i}p_{il}}{\sqrt{\phi^{-2}+|\nabla_{g}f|^{2}}}\nonumber\\
& &+|\nabla_{g}\log\phi|^{2}
+\frac{2f^{l}f^{i}p_{il}f^{j}(\log\phi)_{,j}}{(\phi^{-2}+|\nabla_{g}f|^{2})^{3/2}}
-\frac{\phi^{-2}(\phi^{l}f_{,l})^{2}}{\phi^{-2}+|\nabla_{g}f|^{2}}.\nonumber
\end{eqnarray}
The two expressions (6.20) and (6.21) may now be combined to
obtain
\begin{eqnarray*}
& &\sum_{i}(h_{i4}-K_{i4})^{2}+\sum_{i,j}(h_{ij}-K_{ij})^{2}+2Q(k,\phi,f)\\
&=&-\frac{2\phi^{-1}(\phi^{l}f_{,l})}{\sqrt{\phi^{-2}+|\nabla_{g}f|^{2}}}
(H-\sum_{i}K_{ii})+g^{il}\overline{g}^{jk}p_{ij}p_{lk}\\
& &-2\Delta_{\overline{g}}\log\phi
-|\nabla_{\overline{g}}\log\phi|^{2}
+2\phi^{-1}\Delta_{g}\phi.\nonumber
\end{eqnarray*}
For this last calculation it is necessary to use (6.19) in
addition to
\begin{eqnarray*}
p_{ij}&=&h(X_{i},X_{j})-k_{ij}-f_{,i}f_{,j}k_{44}\\
&=&\frac{1}{\sqrt{\phi^{-2}+|\nabla_{g}f|^{2}}}(f_{;ij}
+(\log\phi)_{,i}f_{,j}+(\log\phi)_{,j}f_{,i})-k_{ij}.
\end{eqnarray*}
The Claim now follows from the generalized Jang equation.
Q.E.D.\medskip

  Theorem 1 is now a consequence of (6.17) as well as Claims 4, 5, 6, and 7.

\begin{center}
\textbf{7.  Appendix B}
\end{center} \setcounter{equation}{0}
\setcounter{section}{7}

   Suppose that the graph map $G(x)=(x,f(x))$ provides an
isometric embedding of $(M,g)$ into the Schwarzschild spacetime
$(\mathbb{SC}^{4},g_{\mathrm{SC}}-\phi^{2}dt^{2})$.  We will show
that $h=K|_{\Sigma}$ implies that the second fundamental form
$\pi$ of the embedding $G(M)\subset\mathbb{SC}^{4}$ is given by
the initial data $k$, where $h$ is the second fundamental form of
the graph $t=f(x)$ (denoted by $\Sigma$) in the warped product
space $(M\times\mathbb{R},g+\phi^{2}dt^{2})$ and $K|_{\Sigma}$ is
the restriction to $\Sigma$ of the extended (by (2.2)) version of
$k$.\par
   Let
\begin{equation*}
X_{i}=\partial_{x^{i}}+f_{,i}\partial_{x^{4}},\text{ }\text{
}\text{ }\text{ }i=1,2,3,
\end{equation*}
be tangent vectors to $\Sigma$.  Then according to (6.14) the
second fundamental form of
$\Sigma\subset(M\times\mathbb{R},g+\phi^{2}dt^{2})$ is given by
\begin{equation*}
h_{ij}:=h(X_{i},X_{j})= \frac{\nabla_{ij}f+(\log\phi)_{,i}f_{,j}
+(\log\phi)_{,j}f_{,i}+g^{lp}\phi\phi_{,l}f_{,p}f_{,i}f_{,j}}{
\sqrt{\phi^{-2}+|\nabla_{g}f|^{2}}}
\end{equation*}
where $\nabla_{ij}$ denotes covariant differentiation with respect
to $g$, and the second fundamental form of
$G(M)\subset\mathbb{SC}^{4}$ is given by
\begin{equation*}
\pi_{ij}:=\pi(X_{i},X_{j})=\frac{\nabla_{ij}^{\mathrm{SC}}
f+(\log\phi)_{,i}f_{,j}+(\log\phi)_{,j}f_{,i}
-g_{\mathrm{SC}}^{lp}\phi\phi_{,l}f_{,p}f_{,i}f_{,j}}
{\sqrt{\phi^{-2}-|\nabla_{g_{\mathrm{SC}}}f|^{2}}}
\end{equation*}
where $\nabla_{ij}^{\mathrm{SC}}$ denotes covariant
differentiation with respect to $g_{\mathrm{SC}}$. Utilizing the
isometry we can write $g_{ij}=(g_{\mathrm{SC}})_{ij}-\phi^{2}
f_{,i}f_{,j}$.  Then direct calculation shows that
\begin{eqnarray*}
& &\nabla_{ij}f+(\log\phi)_{,i}f_{,j}
+(\log\phi)_{,j}f_{,i}+g^{lp}\phi\phi_{,l}f_{,p}f_{,i}f_{,j}\\
&=&\frac{\phi^{-2}}{\phi^{-2}-|\nabla_{g_{\mathrm{SC}}}f|^{2}}(
\nabla_{ij}^{\mathrm{SC}}f+(\log\phi)_{,i}f_{,j}
+(\log\phi)_{,j}f_{,i}).
\end{eqnarray*}
Furthermore
\begin{equation*}
\frac{1}{\sqrt{\phi^{-2}+|\nabla_{g}f|^{2}}}=\phi^{2}
\sqrt{\phi^{-2}-|\nabla_{g_{\mathrm{SC}}}f|^{2}},
\end{equation*}
and so
\begin{equation*}
h_{ij}=\pi_{ij}+\frac{<\phi\nabla_{g}\phi,\nabla_{g}f>_{g}}
{\sqrt{\phi^{-2}+|\nabla_{g}f|^{2}}}f_{,i}f_{,j}.
\end{equation*}
However
\begin{equation*}
(K|_{\Sigma})_{ij}=K(X_{i},X_{j})= k_{ij}+
\frac{<\phi\nabla_{g}\phi,\nabla_{g}f>_{g}}
{\sqrt{\phi^{-2}+|\nabla_{g}f|^{2}}}f_{,i}f_{,j}.
\end{equation*}
Therefore if $h=K|_{\Sigma}$, it follows that $k=\pi$ as desired.

\begin{center}
\textbf{8.  Appendix C}
\end{center} \setcounter{equation}{0}
\setcounter{section}{8}

   Here we calculate the boundary term in (4.4).  Observe that
\begin{equation}
\overline{g}(q,n_{\overline{g}})d\sigma_{\overline{g}}
=\overline{g}^{ij}q_{i}n_{j}\sqrt{g_{11}+\phi^{2}
f_{,r}^{2}}\rho^{2}d\sigma,
\end{equation}
where $(n_{1},n_{2},n_{3})$ is the Euclidean unit outer normal and
$d\sigma$ is the Euclidean area element.  Since the metric
$\overline{g}$ is diagonal and $f=f(r)$ we have $q_{2}=q_{3}=0$.
Next we calculate
\begin{eqnarray}
q_{1}&=&\frac{f^{r}}{\sqrt{\phi^{-2}+|\nabla_{g}f|^{2}}}
(h_{rr}-(K|_{\Sigma})_{rr})\nonumber\\
&=&\frac{\phi g^{11}f_{,r}}{\sqrt{1+\phi^{2}g^{11}f^{2}_{,r}}}
\left(\frac{f_{;rr}+2(\log\phi)_{,r}f_{,r}}{\sqrt{\phi^{-2}+
|\nabla_{g}f|^{2}}}-k_{rr}\right)\nonumber\\
&=&(g_{11}+\phi^{2}f^{2}_{,r})^{-1}\phi f_{,r}(\phi f_{,rr}
-\frac{1}{2}\phi g^{11}g_{11,r}f_{,r}+2\phi_{,r}f_{,r})
-\sqrt{g^{11}}vk_{rr}\nonumber\\
&=&(g_{11}+\phi^{2}f^{2}_{,r})^{-1}\left[\frac{1}{2}g_{11}(\phi^{2}
g^{11}f^{2}_{,r})_{,r}+\phi\phi_{,r}f^{2}_{,r}\right]-\sqrt{g_{11}}vk_{a}\\
&=&g^{11}(1-v^{2})\left[\frac{1}{2}g_{11}\left(\frac{v^{2}}{1-v^{2}}\right)_{,r}
+g_{11}\frac{\phi_{,r}}{\phi}\frac{v^{2}}{1-v^{2}}\right]
-\sqrt{g_{11}}vk_{a}\nonumber\\
&=&\frac{vv_{,r}}{1-v^{2}}+\frac{\phi_{,r}}{\phi}v^{2}
-\sqrt{g_{11}}vK_{a}\nonumber\\
&=&\sqrt{g_{11}}\frac{v}{1-v^{2}}\left[-2\left(\sqrt{g^{11}}\frac{\rho_{,r}}{\rho}
v-k_{b}\right)+(1-v^{2})k_{a}-\sqrt{g^{11}}v\frac{\phi_{,r}}{\phi}(1-v^{2})\right]\nonumber\\
& &+\frac{\phi_{,r}}{\phi}v^{2}-\sqrt{g_{11}}vk_{a}\nonumber\\
&=&-2\sqrt{g_{11}}\frac{v}{1-v^{2}}\left(\sqrt{g^{11}}\frac{\rho_{,r}}{\rho}
v-k_{b}\right)\nonumber,
\end{eqnarray}
where $f_{;rr}$ denotes covariant differentiation with respect to
$g$ and we have used equation (3.2) as well as
$k_{rr}=k(\partial_{r},\partial_{r})=g_{11}k_{a}$. Furthermore
\begin{equation}
\phi=\rho_{,s}=\frac{\sqrt{1-v^{2}}}{\sqrt{g_{11}}}
\rho_{,r},\text{ }\text{ }\text{ }\text{ }\sqrt{g_{11}+\phi^{2}
f_{,r}^{2}}=\frac{\sqrt{g_{11}}}{\sqrt{1-v^{2}}}.
\end{equation}
Thus combining (8.1), (8.2), and (8.3) we have
\begin{equation*}
\phi\overline{g}(q,n_{\overline{g}})d\sigma_{\overline{g}}
=\pm\frac{2\rho_{,r}v}{\sqrt{g_{11}}}
\left(\sqrt{g^{11}}\frac{\rho_{,r}}{\rho}v-k_{b}\right)\rho^{2}d\sigma.
\end{equation*}
\pagebreak

\begin{center}\bigskip
\textbf{References}
\end{center}\bigskip

\noindent 1.\hspace{.07in} H. Bray, \textit{Proof of the
Riemannian Penrose conjecture using the positive mass}
\par\hspace{-.01in} \textit{theorem,} J. Differential Geom.,
$\mathbf{59}$ (2001), 177-267.\medskip

\noindent 2.\hspace{.07in} H. Bray, and M. Khuri, \textit{PDE's
which imply the Penrose conjecture}, preprint,
\par\hspace{-.01in} arXiv:0905.2622, 2009.\medskip

\noindent 3.\hspace{.07in} R. Geroch, \textit{Energy extraction},
Ann. New York Acad. Sci. $\mathbf{224}$ (1973), 108-117.\medskip

\noindent 4.\hspace{.07in} S. Hayward, \textit{Gravitational
energy in spherical symmetry}, Phys. Rev. D, $\mathbf{53}$
\par\hspace{-.01in} (1996), 1938-1949.\medskip

\noindent 5.\hspace{.07in} S. Hayward, \textit{Inequalities
relating area, energy, surface gravity, and charge of black}
\par\hspace{-.01in} \textit{holes}, Phys. Rev. Lett., $\mathbf{81}$
(1998), 4557-4559.\medskip

\noindent 6.\hspace{.07in} G. Huisken, and T. Ilmanen, \textit{The
inverse mean curvature flow and the Rieman-}
\par\hspace{-.01in} \textit{nian Penrose inequality}, J. Differential Geom.,
$\mathbf{59}$ (2001), 353-437.\medskip

\noindent 7.\hspace{.06in} M. Iriondo, E. Malec, and N. \'{O}
Murchadha, \textit{Constant mean curvature slices and}
\par\hspace{-.01in} \textit{trapped surfaces in asymptotically flat
spherical spacetimes}, Phys. Rev. D, $\mathbf{54}$
\par\hspace{-.01in} (1996), 4792-4798.\medskip

\noindent 8.\hspace{.07in} P.-S. Jang, \textit{On the positivity
of energy in General Relativity}, J. Math. Phys., $\mathbf{19}$
\par\hspace{-.01in} (1978), 1152-1155.\medskip

\noindent 9.\hspace{.07in} M. Khuri, \textit{A Penrose-like
inequality for general initial data sets}, Commun. Math.
\par\hspace{-.01in} Phys., $\mathbf{290}$ (2009), 779-788.\medskip

\noindent 10.  E. Malec, and N. \'{O} Murchadha, \textit{Trapped
surfaces and the Penrose inequality in} \par\hspace{-.01in}
\textit{spherically symmetric geometries}, Phys. Rev. D,
$\mathbf{49}$ (1994), 6931-6934.\medskip

\noindent 11.  E. Malec, and N. \'{O} Murchadha, \textit{The Jang
equation, apparent horizons, and the}
\par\hspace{-.01in} \textit{Penrose inequality}, Class. Q. Grav.,
$\mathbf{21}$ (2004), 5777-5787.\medskip

\noindent 12.  J. Metzger, \textit{Blowup of Jang's equation at
outermost marginally trapped surfaces},
\par\hspace{-.01in} preprint, arXiv:0711.4753, 2008.\medskip

\noindent 13.  R. Schoen, and S.-T. Yau, \textit{Proof of the
positive mass theorem II}, Commun. Math.
\par\hspace{-.01in} Phys., $\mathbf{79}$ (1981), 231-260.

\bigskip\bigskip\footnotesize

\noindent\textsc{Department of Mathematics, Duke University, Box
90320, Durham, NC 27708}\par

\noindent\textit{E-mail address}: \verb"bray@math.duke.edu"
\bigskip

\noindent\textsc{Department of Mathematics, Stony Brook
University, Stony Brook, NY 11794}\par

\noindent\textit{E-mail address}: \verb"khuri@math.sunysb.edu"

\end{document}